\documentclass[reqno]{amsart}

\usepackage{amsmath,amsthm,amssymb,amsfonts,mathrsfs,amscd}
\usepackage{multirow}
\usepackage{graphicx}

\usepackage{bm}
\usepackage{mathrsfs}
\usepackage{amsmath}
\usepackage{graphicx}
\usepackage{amssymb}
\usepackage{dsfont}
\usepackage{tikz}
\usepackage{graphicx}
\usepackage{tikz-3dplot}
\usepackage{subfig}
\usepackage{bigstrut}

\newtheorem{remark}{Remark}[section]

\newtheorem{lemma}{Lemma}[section]
\newtheorem{theorem}{Theorem}[section]

\newcommand{\stl}
{\mathrel{\raise2pt\hbox{${\mathop<\limits_{\raise1pt\hbox
{\mbox{$\sim$}}}}$}}}

\newcommand{\ste}
{\mathrel{\raise2pt\hbox{${\mathop=\limits_{\raise1pt\hbox{\mbox{$\sim$}}}}$}}}

\newcommand{\stg}
{\mathrel{\raise2pt\hbox{${\mathop>\limits_{\raise1pt\hbox{\mbox{$\sim$}}}}$}}}

\makeatletter
  \newcommand\figcaption{\def\@captype{figure}\caption}
  \newcommand\tabcaption{\def\@captype{table}\caption}
\makeatletter

\def \no{\noindent}
\def\su{\sum\limits^{N_0}_{i=1}}

\def\be{\begin{equation}}
\def\bn{\begin{eqnarray}{}}
\def\ee{\end{equation}}
\def\en{\end{eqnarray}}
\def\bq{\begin{eqnarray*}{}}
\def\eq{\end{eqnarray*}}
\def\l{\label}
\def\ri{\rightarrow}
\def\o{\omega}
\def\O{\Omega}
\def\la{\lambda}
\def\al{\alpha}

\def\T{{\mathcal T}_h}

\def\G{\Gamma}
\def\g{\gamma}
\def\p{\partial}

\def\12{{1\over 2}}

\def\a{{\mathcal A}}
\def\q{\quad}

\def\V0{\bar{V}_h}

\def\ff{ {\mbox{\sc f}}}
\def\gg{ {\mbox{\sc g}}}
\def\se{{\mbox{\sc e}}}
\def\vv{\mbox{\sc v}}
\def\o{\Omega}
\def\c{{\bf curl}}
\def\u{{\bf u}}
\def\v{{\bf v}}

\def\x{{\bf x}}
\def\ti{\times}

\def\g{\Gamma}

\def\w{{\bf w}}
\def\n{{\bf n}}
\def\f{{\bf f}}
\def\0{{\bf 0}}
\def\a{{\bf a}}
\def\b{{\bf b}}
\def\t{{\bf t}}
\def\r{{\bf r}}

\def\12{{1\over 2}}

\def\su{\sum\limits}

\def\rig{\rightarrow}

\def\chi{{\mathcal X}}

\def\A{{\mathcal A}}

\begin{document}

\title{A substructuring preconditioner with vertex-related interface solvers for
elliptic-type equations in three dimensions}

\author{Qiya Hu}
\author{Shaoliang Hu}

\thanks{1. LSEC, ICMSEC, Academy of Mathematics and Systems Science, Chinese Academy of Sciences, Beijing
100190, China; 2. School of Mathematical Sciences, University of Chinese Academy of Sciences, Beijing 100049,
China (hqy@lsec.cc.ac.cn \and hushaoliang@lsec.cc.ac.cn).
This work was funded by Natural Science Foundation of China G11571352.}

\maketitle
\begin{abstract}
In this paper we propose a variant of the substructuring preconditioner for solving
three-dimensional elliptic-type equations with strongly discontinuous coefficients.
In the proposed preconditioner, we use the simplest coarse solver associated with the finite element space induced by the coarse
partition, and construct inexact interface solvers based on overlapping domain decomposition with small overlaps. This new preconditioner has
an important merit: its construction and efficiency
do not depend on the concrete form of the considered elliptic-type equations.
We apply the proposed preconditioner to solve the linear elasticity problems and Maxwell's equations in three dimensions.
Numerical results show that the convergence rate of PCG method
with the preconditioner is nearly optimal, and also robust with respect to
the (possibly large) jumps of the coefficients in the considered equations.
\end{abstract}



\pagestyle{myheadings}
\thispagestyle{plain}
\markboth{Substructuring preconditioners  for elliptic-type equations}{QIYA HU AND SHAOLIANG HU}

\section{Introduction}
There are many works to study domain decomposition methods (DDMs) for solving the systems generated by finite element discretization of
elliptic-type partial differential equations (\cite{BrambleP1989}-\cite{CaiS1999},\cite{ChanZ1994}-\cite{Gander2012}, \cite{Hu2004}-\cite{Mandel2005}, \cite{Smith1992}-\cite{Toselli2006}, \cite{VeigaCPS2012}-\cite{VeigaPSWZ2014},\cite{re25,XuZ1998} and the references therein). It is known that, for three dimensional problems, non-overlapping domain decomposition methods (DDMs) are more difficult to construct and implement than overlapping DDMs although the non-overlapping DDMs have some advantages over the overlapping DDMs in the treatment of jump coefficients.
In fact, the construction of non-overlapping DDMs heavily depends on the considered models in three dimensiona. For example, non-overlapping DDMs for positive definite Maxwell's equations are essentially different
from the usual elliptic equation (comparing  \cite{DohrmannW2015,HuSZ2011,HuZ2004,Toselli2006}). This drawbacks restrict applications of the non-overlapping DDMs.

A key ingredient in the construction of non-overlapping domain decomposition methods is the choice of a suitable coarse subspace. There are two main ways to construct
coarse subspaces in the existing works: (i) use some degrees of freedom on the {\it joint-set} (BPS method, FETI-DP method, BDDC method); (ii) use the local kernel spaces of the considered differential operator
(Neumann-Meumann method, BDD method, FETI method). But, there are some drawbacks in the two ways. For the first way, the choice of degrees of freedom heavily depends on the considered model, for example, one can choose the degrees of freedom on the vertices, the averages on the edges or faces for the three dimensional Laplace equation, but this choice is not practical for the three dimensional Maxwell's equations
(see \cite{DohrmannW2015}, \cite{HuZ2003}, \cite{HuZ2004} and \cite{Toselli2006}). For the second way,
the coarse space may be very large. For example, the kernel space of the $curl$ operator is just the gradient of the nodal
finite element space, so (by simple calculation) the dimensions of such coarse space for Maxwell's equations are greater than $1/7$ of that of the original solution space.

Another possible choice of coarse subspace is the finite element space induced by the coarse partition. This coarse subspace
was first considered in \cite{DryjaSW1994} for elliptic equations, and then was investigated in \cite{XuZ1998}.
It is clear that this coarse subspace possesses the simplest structure and almost the smallest degrees of freedom among the coarse subspaces considered in the existing non-overlapping DDMs
for three-dimensional problems. However, such coarse subspace
was regarded as a unapplicable coarse subspace for long time, since the condition number of the resulting preconditioned system is not nearly optimal for three dimensional problems
with large jump coefficients. For elliptic equations in three dimensions, substructuring preconditioners with such coarse subspace was studied again in \cite{HuS2010}.
Based on the framework developed in \cite{XuZ2008}, it was shown that the PCG method for solving the resulting preconditioned system has the
nearly stable convergence even for the case with large jump coefficients. In \cite{HuSZ2011}, this kind of coarse solver was also applied to the construction of substructuring preconditioner
 for Maxwell's equations in three dimensions. Unfortunately, for this choice of coarse
subspace, we need to design local interface solvers on coarse edges, whose definitions still depend on the considered models (comparing \cite{HuSZ2011} with \cite{HuS2010}).

In the present paper, we try to construct relatively united substructuring preconditioner for elliptic-type equations, such that it is cheap, easy to implement and has fast convergence.
As usual, we decompose the considered domain into the union of some non-overlapping subdomains, which constitute a coarse partition of the domain.
In the proposed preconditioner, we use the simplest coarse space induced by the coarse
partition (see \cite{DryjaSW1994}, \cite{HuS2010}, \cite{HuSZ2011} and \cite{XuZ1998}). The main goal of this paper is to design unified and practical local interface solvers.

For each internal cross-point, we introduce an auxiliary subdomain that contains the internal cross-point as its ``center" and has almost the same size with the original subdomains. Associated
with each auxiliary subdomain, we define a local interface solver such that the solution of the local interface problem is discrete harmonic in the intersection of the auxiliary subdomain with every original subdomain
adjoining to it. Notice that each intersection is only a part of some original subdomain, so the local interface solver corresponds to ``inexact" harmonic extensions. Such a local interface solver is implemented by
solving a Dirichlet problem (residual equation), which is defined on the natural restriction space of the original finite element space on the auxiliary subdomain. It is clear that each local interface solver
has almost the same cost with an original subdomain solver. We point out that the proposed local interface solvers are different from the existing local interface solvers defined in the vertex space method \cite{Smith1992} or the interface overlapping additive Schwarz  \cite{XuZ1998}, where exact harmonic extensions are required.

In order to further reduce the cost of the local interface solvers described above,
we present approximate local interface solvers based on a coarsening technique. In the step for
solving a local interface problem, we are interested only in the degrees of freedom on
the local interface, instead of the degrees of freedom in the interiors of subdomains.
Intuitively, the accuracy of the degrees of freedom on the local interface are not
sensitive to the grids far from the local interface. Based on this observation, we construct
auxiliary non-uniform grids in each subdomain adjoining the considered local
interface such that the auxiliary grids coincide with the original fine grids on the local
interface but gradually become coarser when nodes are far from the local interface.The desired
approximate local interface solver is implemented by solving a Dirichlet problem on
the finite element space defined by the auxiliary grids, which have much smaller nodes than the original fine grids.

The constructions of the coarse solver and the proposed local interface solvers
do not depend on the considered models, and the resulting substructuring preconditioner is cheap and easy to implement.
As pointed out in \cite{DohrmannW2015}, the design of an efficient substructuring preconditioner for three dimensional Maxwell's equations poses quite significant
challenges. A few existing preconditioners on this topic are either expensive or difficult to implement. We will apply the proposed substructuring preconditioner to solve the linear elasticity problems
and Maxwell's equations in three dimensions. Numerical results show that the preconditioner is robust uniformly for the two kinds of equations even if the coefficients have large jumps.
We also consider possible extension of the preconditioner to the case with irregular subdomains.

The outline of the paper is as follows. In Section 2, we give the variational formula of general elliptic-type equations and introduce a partition based on domain
decomposition.  In Section 3, we describe local interface solvers associated with vertex-related subdomains and define the resulting substructuring preconditioner for the general elliptic system.
In Section 4, we present an analysis of convergence of the preconditioner for linear elasticity problems.
In section 5, we will report the some numerical results for the linear elasticity problems and Maxwell's
equations.

\section{Elliptic-type equations and domain decomposition}
In this section, we describe the considered problems.
\subsection{Elliptic-type equations}
Let $\Omega$ be a bounded and connected Lipschizt domain in $\mathds{R}^3$. For convenience, we just consider the weak form of elliptic-type equations. Let $V(\Omega)$ denote a Hilbert
space with the scalar product $(\cdot,\cdot)_V$, and $||\cdot||_V$ be the induced norm.
We introduce a real bilinear form $\A(\cdot,\cdot): V(\Omega) \times V(\Omega) \rightarrow R$. We assume that $\A(\cdot,\cdot)$ is symmetric,
\[ \A(\u,\v) = \A(\v,\u), \quad \u,\v \in V(\Omega) \]
continuous,
\[ |\A(\u,\v)| \leq c_1||\u||_V ||\v||_V, \quad \u,\v \in V(\Omega) \quad c_1 >0 \]
and coercive
\[ \A(\u,\u) \geq c_2 ||\u||^2_V, \quad \u,\v\in V(\Omega), \quad c_2 > 0. \]

Giving a linear functional $ \bm{F} \in V'(\Omega)$, we consider the following problem:
\begin{equation}
\begin{cases}
 Find \q \u \in V(\o) \quad .st. \\
\A(\u, \v)=(\bm{F},\v), \quad \forall\v\,\in V(\Omega),
 \end{cases}
 \label{eq:2.1}
\end{equation}
where $(\cdot, \cdot)$ denotes the duality pairing between $V'(\Omega)$ and $V(\Omega)$.

\subsection{Domain decomposition and discretization}
For convenience, we assume that $\Omega$ is a polyhedra. For a number $d\in
(0,~1)$, let $\Omega$ be decomposed into the union of
non-overlapping tetrahedra (or hexahedra) $\{\O_k\}$ with the size
$d$. Then, we get a non-overlapping domain decomposition for $\O$:
$\bar{\O}=\bigcup\limits_{k=1}^{N} \bar{\O}_k$. Assume that
$\O_i\cap\O_j=\emptyset$ when $i\not=j$; if $i\not=j$ and
$\partial\O_i\cap\partial\O_j\not=\emptyset$, then
$\partial\O_i\cap\partial\O_j$ is a common face of $\O_i$ and
 $\O_j$, or a common edge of $\O_i$ and
 $\O_j$, or a common vertex of $\O_i$ and
 $\O_j$. It is clear that the subdomains $\O_1,\cdots,\O_N$
constitute a {\it coarse} partition $\mathcal{T}_d$ of $\O$. If
$\partial\O_i\cap\partial\O_j$ is just a common face of $\O_i$ and
 $\O_j$, then set $\G_{ij}=\partial\O_i\cap\partial\O_j$. Define $\G=\cup\G_{ij}$.
 By $\Gamma_k$ we denote the intersection of $\Gamma$ with the boundary of the subdomain
$\o_k$. So we have $\Gamma_k=\p\o_k$ if $\o_k$ is an interior
subdomain of $\o$.

With each subdomain $\O_k$ we associate a regular triangulation made
of tetrahedral elements (or hexahedral elements). We require that
the triangulations in the subdomains match on the interfaces between
subdomains, and so they constitute a triangulation $\T$ on the
domain $\O$, which we assume is quasi-uniform. We denote by $h$ the mesh
 size of $\T$, i.e., $h$ denotes the maximum diameter of tetrahedra in
 the mesh $\T$.

For an element $K\in {\mathcal T}_h$, let $R(K)$ denote a set of basis functions on the element $K$. The definition
of $R(K)$ depends on the considered models, and will be given in Section 4.2 and Subsection 5.2. Define the finite element space
$$ V_{h}(\o)=\Big\{\v\in V({\Omega)}: ~ \v|_K\in R(K), ~\forall K\in \mathcal{T}_h
 \Big\},
$$
Consider the discrete problem of (\ref{eq:2.1}):
{\it Find $\bm{u}_h\in V_{h}(\O)$ such that}
\begin{equation}
\A(\bm{u}_h, \bm{v})=(\bm{F}, \bm{v}),
\quad\forall \bm{v}\in V_{h}(\O). \label{eq:2.2}
\end{equation}
This is the whole problem we need to solve in this paper.

For convenience, we define the discrete operator $A: V_h(\Omega) \rightarrow V_h(\Omega)$ as
\[ (A\u, \v) = \A(\u,\v), \quad \u, \v \in V_h(\Omega). \]
Then (\ref{eq:2.2}) can be written in the operator form
\begin{equation}
Au_h=\bm{f}.
\label{eq:2.3}
\end{equation}
By the assumptions on ${\mathcal A}(\cdot,\cdot)$, the operator is symmetric and positive definite. Thus the above equation can be iteratively solved by CG method.
In the rest of this paper, we will construct a preconditioner for the operator $A$.

Before constructing the desired preconditioner, we first introduce some useful sets and subspaces.

$\mathcal{N}_h$: the set of all nodes generated by the {\it fine} partition $\mathcal{T}_h$;

$\mathcal{E}_h$: the set of all {\it fine} edges generated by the partition $\mathcal{T}_{h}$;

$\mathcal{F}_h$: the set of all {\it fine} faces generated by the partition
$\mathcal{T}_h$;

$\mathcal{N}_d$: the set of all nodes generated by the {\it coarse} partition $\mathcal{N}_d$.




In most applications, the degrees of freedom of $\v\in V_{h}(\o)$ are defined at the nodes in $\mathcal{N}_h$ (the nodal elements), or on the edges in $\mathcal{E}_h$ (Nedelec edge elements), or on the
faces in $\mathcal{F}_h$ (Raviart-Thomas face elements).  Throughout this paper, for a subset $\ff$ that is the union of faces in $\mathcal{F}_h$, the term
``the~degrees~of~freedom~of $\v$ vanish on~$\ff$" means that ``$\v$ has the zero~degrees~of~freedom~at~the nodes, or fine edges, or fine faces of $\ff$".

Let $G\subset \Omega$ be a subdomain that is the union of some elements in ${\mathcal T}_h$. Define
$$ V^0_h(G)=\{v\in V_{h}(\o):~~\mbox{the~degrees~of~freedom~of}~\v~\mbox{vanish~on}~\partial G\}.$$
For example, when $G=\Omega_k$ the space $V^0_h(\Omega_k)$ is just the subdomain space in the traditional substructuring methods.

For the construction of the desired preconditioner, we will use the simplest coarse space $V_d(\Omega)$, which is defined as the finite element space
associated with the {\it coarse} partition $\mathcal{T}_d$ (see \cite{DryjaSW1994}, \cite{HuS2010}, \cite{HuSZ2011} and \cite{XuZ1998}). It is clear that
$V_{d}(\o)\subset V_{h}(\o)$.

\section{A preconditioner with vertex-related local interface solvers}
\setcounter{equation}{0}
This section is devoted
to describing the desired preconditioner, in which local interface solvers are defined in vertex-related subspaces.

\subsection{Space decomposition}
\begin{figure}[ht]
\begin{center}
\includegraphics[width =7cm]{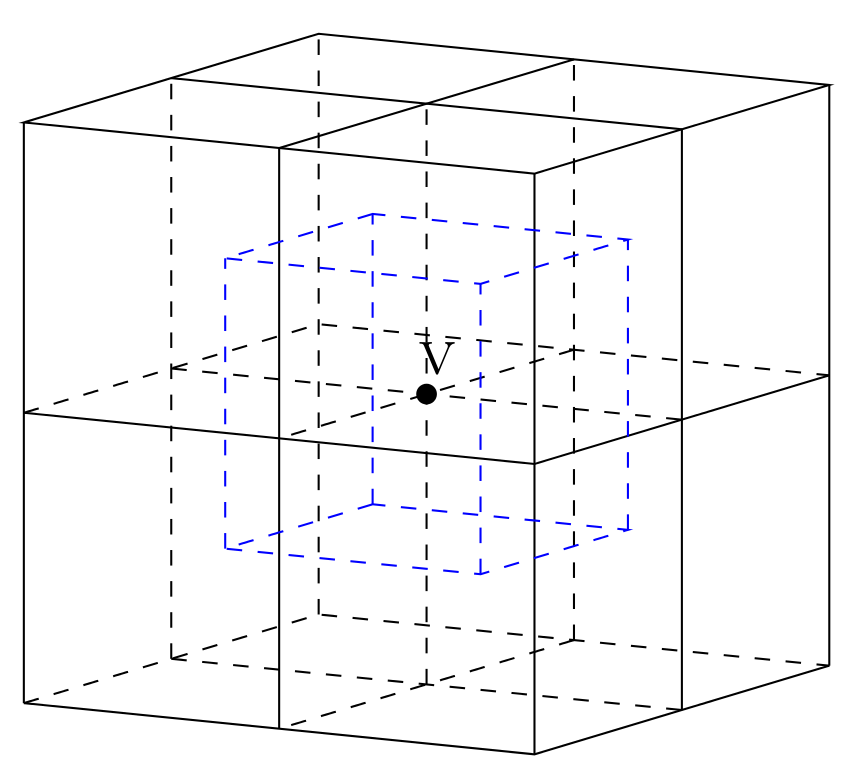}   
\caption{The auxiliary subdomain $\Omega_{\vv}^{half}$ (the blue cube) associated with the vertex $\vv$. \label{map1}}
\end{center}
\end{figure}

For each $\vv \in \mathcal{N}_d$, we construct an open region $\Omega_{\vv}^{half}$, whose ``center" is $\vv$ and size is about $d$, see Figure \ref{map1}. When $\vv\in\partial\Omega$,
the auxiliary subdomain
$\Omega_{\vv}^{half}$ is chosen as the part in $\Omega$. We assume that: (i) each
subdomain $\Omega_{\vv}^{half}$ is just the union of some elements in ${\mathcal T}_h$; (ii) the union of all the interface subdomains $\Omega_{\vv}^{half}\cap\Gamma$ is an open cover of
$\Gamma$. Then all the interface subdomains $\Omega_{\vv}^{half}\cap\Gamma$ constitute an overlapping domain decomposition of $\Gamma$ (with small overlap), where the overlap can be one element
layer only. Here we do not require that all the subdomains $\Omega_{\vv}^{half}$ constitute an overlapping domain decomposition of the original domain $\Omega$, so we can choose slightly smaller
subdomains $\Omega_{\vv}^{half}$ in applications.

In order to define a space decomposition of $V_h(\Omega)$ in an exact manner, we need to introduce more notations.

For $\vv\in {\mathcal N}_d$, set
$$ \Lambda_{\vv}=\{k:~~\mbox{the~polyhedran}~\Omega_k~\mbox{contains}~\vv~\mbox{as~its~ vertex}\} $$
and define
$$ \Gamma_{\vv}^{half}=\Omega_{\vv}^{half}\cap\Gamma~~\mbox{and}~~\Omega_{\vv}=\bigcup_{k\in\Lambda_{\vv}}\Omega_k.$$
Let $V_h(\Gamma)$ denote the interface space, which consists of the natural traces of all the functions in $V_h(\Omega)$. Define the vertex-related local interface space
$$ V_h^0(\Gamma_{\vv}^{half})=\{\phi\in V_h(\Gamma):~~supp~\phi\subset \Gamma_{\vv}^{half}\}. $$
Since all the vertex-related local interfaces $\Gamma_{\vv}^{half}$ constitute an open cover of the interface $\Gamma$, we have the space decomposition
\begin{equation}
V_h(\Gamma)=\bigcup_{\vv\in{\mathcal N}_d}V^0_h(\Gamma^{half}_{\vv}).\label{3.new00}
\end{equation}

As usual, let $V_h^{\bot}(\Omega)$ denote the space consisting of all the finite element functions that are discrete harmonic in each $\Omega_k$, namely
$$ V^{\bot}_h(\Omega)=\{\v\in V_h(\o):~~{\mathcal A}(\v,\w)=0,~\forall\w\in V_h^{0}(\o_{k})~~{for~each}~~\Omega_k\}.$$
For each $\vv\in {\mathcal N}_d$, define vertex-related local harmonic space
$$ V^{\bot}_h(\Omega_{\vv})=\{\v\in V^{\bot}_h(\o):~~\mbox{the~trace~of}~\v~\mbox{belongs~to}~V_h^0(\Gamma_{\vv}^{half})\}\subset V_h^0(\Omega_{\vv}).$$
In other words, $V^{\bot}_h(\Omega_{\vv})$ is just the space consisting of the discrete harmonic extensions of the functions in $V_h^0(\Gamma_{\vv}^{half})$.

It is clear that
$$ V^{\bot}_h(\Omega)=\bigcup_{\vv\in{\mathcal N}_d}V^{\bot}_h(\Omega_{\vv}). $$
Thus the space $V_h(\o)$ admits the space decomposition
\begin{equation}
 V_h(\o)=V_d(\o)+\su_{k=1}^NV_h^0(\o_k)+\su_{\vv\in\mathcal{N}_d}V^{\bot}_h(\Omega_{\vv}).\label{eq:3.1}
\end{equation}

\subsection{Preconditioner}
In this subsection we define solvers on the subspaces $V_d(\o)$,
$V_h^0(\o_k)$ and $V^{\bot}_h(\o_{\vv})$.

As usual, we use $A_d: V_d(\o)\rig V_d(\o)$ and $A_k: V^0_h(\o_k)\rig V_h^0(\o_k)$ to denote the restriction of $A$ on $V_d(\o)$ and $V_h^0(\o_k)$ respectively, i.e., they satisfy
$$ (A_d\v_d,\w_d)=\A(\v_d,\w_d),~~~\v_d\in V_d(\o),~~\forall \w\in V_d(\o) $$
and
$$ (A_k\v, {\w})_{\o_k}=(A\v, {\w})=\A(\v,{\w}),~\v\in V^0_h(\o_k),~\forall {\w}\in V_h^0(\o_k).$$

In the following we define an ``inexact" solver on $V^{\bot}_h(\o_{\vv})$. To this end, we introduce a modification of $V^{\bot}_h(\o_{\vv})$.  Let $k\in\Lambda_{\vv}$, and use $\o_{\vv,k}^{half}$ to denote
the intersection of $\o_{\vv}^{half}$ with $\Omega_k$. For each $\o_{\vv}^{half}$, define the ``inexact" harmonic space
$$ V^{\bot}_h(\o_{\vv}^{half})=\{\v\in V^0_h(\o_{\vv}^{half}):~~{\mathcal A}(\v,\w)=0,~\forall\w\in V_h^{0}
(\o_{\vv,k}^{half})~~\mbox{with}~~k\in\Lambda_{\vv}\}.$$
Notice that the functions in $V^{\bot}_h(\o_{\vv}^{half})$ have the support set $\o_{\vv}^{half}$ and are harmonic only in the subdomain $\o_{\vv,k}^{half}$ of $\Omega_k$ (for any $k\in\Lambda_{\vv}$).

For a function $\v\in V^{\bot}_h(\o_{\vv})$, define $\v^{half}\in V^{\bot}_h(\o_{\vv}^{half})$ such that $\v^{half}=\v$ on $\Gamma_{\vv}^{half}$.
For each $\vv \in \mathcal {N}_d$, let $B_{\vv}: V^{\bot}_h(\o_{\vv})\rig V^{\bot}_h(\o_{\vv})$ be the symmetric and
positive definite operator defined by
$$
(B_{\vv}\v, \w)=\A(\v^{half}, \w^{half}),\q\v\in V^{\bot}_h(\o_{\vv}),~~\forall\w\in V^{\bot}_h(\o_{\vv}). $$
Since the basis functions in
$V^{\bot}_h(\o_{\vv})$ are not known, the action of $B^{-1}_{\vv}$ needs to be implemented by solving a residual equation defined in $V^0_h(\o_{\vv}^{half})$
(see {\bf Algorithm 3.1} given later).

Let $Q_d: V_h(\o)\rig V_d(\o)$, $Q_k:V_h(\o)\rig V^0_h(\o_k)$
and $Q_{\vv}: V_h(\o)\rig V^{\bot}_h(\o_{\vv})$ be the standard $L^2$-projectors.
Then the preconditioner for $A$ is defined as follows:
\be
B^{-1}=A_d^{-1}Q_d+\su_{k=1}^NA^{-1}_kQ_k+\su_{\vv\in \mathcal {N}_d}B_{\vv}^{-1}Q_{\vv}
\l{precon1} \ee

The action of the preconditioner $B^{-1}$, which is needed in each iteration step of
PCG method, can be described by the following algorithm.\\
{\bf Algorithm 3.1}. For ${\bf g}\in V_h(\o)$, we can compute
$\u=B^{-1}{\bf g}$
in four steps.

Step 1. Solve the system for $\u_d\in V_d(\o)$:
$$ (A_d\u_d, \v_d)=({\bf g},\v_d),~~~\forall\v_d\in V_d(\o); $$

Step 2. Solve the following system for  $\u_k\in V_h^0(\o_k)$ ($k=1,\cdots,N$) in parallel:
$$ (A_k\u_k, \v)=({\bf g},\v),~~~\forall\v\in V^0_h(\o_k),~~k=1,\cdots,N;
$$

Step 3. Solve the following system for $\u_{\vv}\in V_h^0(\o_{\vv}^{half})$ ($\vv \in \mathcal{N}_d$) in parallel:
$$ (B_{\vv}\u_{\vv}, \v)=({\bf g},\v)-\su_{k\in\Lambda_{\vv}}(A_k\u_k,\v)_{\o_k},~~~\forall\v\in
V^0_h(\o_{\vv}^{half});$$

Step 4. With the trace
$\Phi_h=\bm{\gamma}(\su_{\vv \in \mathcal{N}_d}\u_{\vv})|_{\g}$,
compute the $A$-harmonic extension of $\Phi_h$ on each $\o_k$ to obtain
$\u^{\bot}\in V^{\bot}_h(\o)$. This leads to
$$ \u=\u_d+\su_{k=1}^N\u_k+\u^{\bot}. $$

\begin{remark} Notice that all the subproblems in {\bf Algorithm 3.1} correspond to the same bilinear form as ${\mathcal A}(\cdot,\cdot)$ (but with different
finite element spaces). We point out that each vertex-related space $V^0_h(\o_{\vv}^{half})$ has almost the same
degrees of freedom with an original subdomain space $V^0_h(\o_k)$. Moreover, the local problem in Step 4 has the same stiffness matrix with that in Step 2
(with different right hands only). Thus the implementation of Step 4 is very cheap by using LU decomposition made in Step 2 for each local stiffness matrix (if Step 2 is implemented in the direct method).
All this shows that {\bf Algorithm 3.1} is easy and cheap to implement.
\end{remark}

\subsection{An approximation of the vertex-related solver $B_{\vv}$}
In order to reduce the cost for the implementation of Step 3 in {\bf Algorithm 3.1}, we would like to replace each space $V_h^0(\o_{\vv}^{half})$ by a smaller space.
In Step 4 we need only to use the values of $\u_{\vv}$ on $\Gamma_{\vv}^{half}$, so we only hope to get a rough approximation of $\u_{\vv}|_{\Gamma_{\vv}^{half}}$
but do not care for the accuracy of $\u_{\vv}$ at the nodes in the interior of $\Omega^{half}_{\vv}\backslash\Gamma_{\vv}^{half}$.
Intuitively, the accuracy of an approximation for $\u_{\vv}$ mainly depends on the grids nearing $\Gamma_{\vv}^{half}$ and is not sensitive to the grids far from $\Gamma_{\vv}^{half}$.
Based on this observation, we can construct an auxiliary non-uniform partition $\tilde{\mathcal T}^{\vv}_{\tilde{h}}$ on $\Omega^{half}_{\vv}$, for which the original
fine grid on $\Gamma_{\vv}^{half}$ are kept and the grid in the interior of $\Omega^{half}_{\vv}\backslash\Gamma_{\vv}^{half}$
gradually becomes coarser when nodes are far from $\Gamma_{\vv}^{half}$. The auxiliary partition $\tilde{\mathcal T}^{\vv}_{\tilde{h}}$ can be easily generated by the existing software \cite{gmsh}.
With this auxiliary partition $\tilde{\mathcal T}^{\vv}_{\tilde{h}}$,
we can define a new linear finite element space $\tilde{V}_{\tilde{h}}^0(\o_{\vv}^{half})$, each function in which vanishes outside $\Omega_{\vv}^{half}$.
Now we replace the space $V_h^0(\o_{\vv}^{half})$ in Step 3 by $\tilde{V}_{\tilde{h}}^0(\o_{\vv}^{half})$ and
we get a variant Step $3'$ of Step 3. The resulting solver, which is denoted by $\tilde{B}_{\vv}$, can be viewed as an approximation of $B_{\vv}$.
Since the dimension of $\tilde{V}_{\tilde{h}}^0(\o_{\vv}^{half})$ is much smaller than that of $V_h^0(\o_{\vv}^{half})$, the implementation of
Step $3'$ (i.e., the action of $\tilde{B}^{-1}_{\vv}$) is much cheaper than that of Step 3 (i.e., the action of $B^{-1}_{\vv}$). For convenience, we use $\tilde{B}$ to denote the
preconditioner defined by (\ref{precon1}) with $B_{\vv}$ being replaced by $\tilde{B}_{\vv}$.

\subsection{An extension to the case with irregular subdomains}
\begin{figure}[hbt]
\centering
\begin{tabular}{c}
\includegraphics[width =8cm]{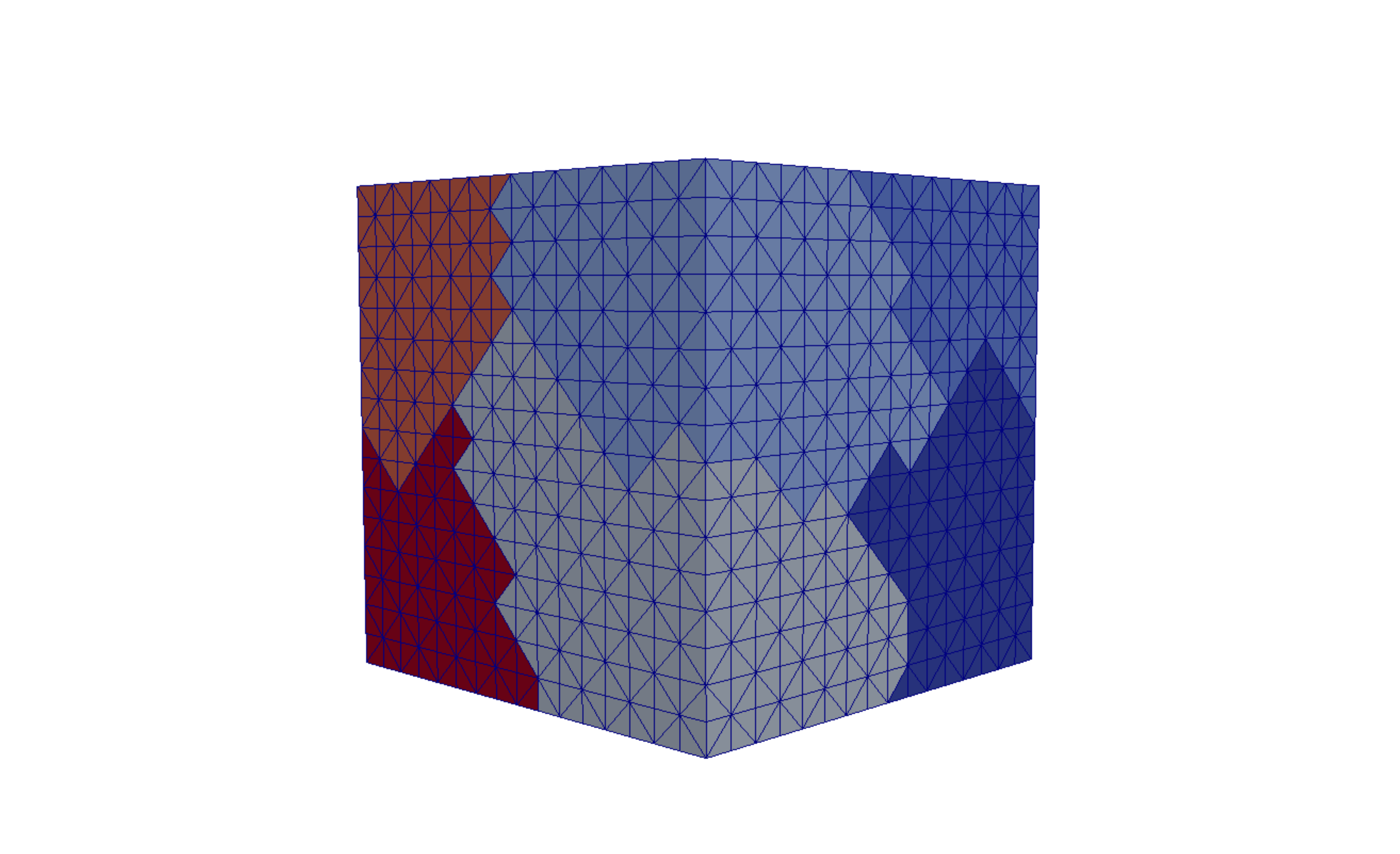}
\end{tabular}
\caption{irregular partition for $\Omega$}
\label{fig:cube_dom}
\end{figure}

In this subsection, we consider the case that the subdomains $\{\Omega_k\}$ in the previous sections
are irregular, i.e., some subdomains $\Omega_k$ are not polyhedrons with finite faces, see Figure \ref{fig:cube_dom}. This situation appears when subdomains
$\{\Omega_k\}$ are automatically generated by the software Metis \cite{Metis} for given fine meshes.

Notice that the subdomains generated by this software may be geometrically non-conform, i.e., the common part of two neighboring subdomains is not
a complete face of one of the two subdomains. In this situation, the information on substructuring vertices, edges and faces are in general unknown and can not be directly obtained.
This means that, before implementing the proposed preconditioner for this case, one needs
to get information on coarse vertices (information on coarse edges and coarse faces is also needed for the BDDC method or the FETI-DP method).
In order to extend the proposed preconditioner to the case with irregular subdomains, we recall a precise definition of substructuring vertices
(refer to \cite{Mandel2003}) and give a variant of the coarse subspace $V_d(\Omega)$.

$\bullet$ Substructuring vertices

 We first extend the definition in \cite{Mandel2003}. Let $\Gamma$ denote the interface, which is the union of all the common parts of two neighboring subdomains,
 and let $\mathcal{N}_{\Gamma}$ denote the set of nodes on $\Gamma$. We need to decompose the set $\mathcal{N}_{\Gamma}$ into the union of disjoint equivalence classes $\{\mathcal{N}_l\}$.
 For this purpose, we define an index set of subdomains for each $p\in\mathcal{N}_{\Gamma}$ by
 $$ \Lambda(p)=\{k:~~p\in \partial\Omega_k\}.$$
Two nodes $p,q\in\mathcal{N}_{\Gamma}$ belong to a same equivalence class $\mathcal{N}_l$ if and only if $\Lambda(p)= \Lambda(q)$.
 For example, all the nodes on the common part of two neighboring subdomains constitute a class $\mathcal{N}_l$.
In particular, some class $\mathcal{N}_l$ may contain only one node on $\Gamma$. For convenience, let ${\mathcal N}_d^{\ast}$ denote the union of all the single point sets $\mathcal{N}_l$,
each of which contains only one node in $\mathcal{N}_{\Gamma}$.

For the standard domain decomposition with polyhedral subdomains, the set ${\mathcal N}_d^{\ast}$ is just the subdomain vertex set used to define coarse subspace.
However, the geometric structure of subdomains generated by Metis is very complicated, which
makes ${\mathcal N}_d^{\ast}$ contain too many vertices and some vertices in ${\mathcal N}_d^{\ast}$ to be very close. In order
to reduce the number of vertices and
construct a practical coarse space that contains small degree of freedoms, we need to choose partial vertices from
${\mathcal N}_d^{\ast}$ to get the desired vertex set ${\mathcal N}_d\subset {\mathcal N}_d^{\ast}$. The main idea is to remove the vertices that have very short distances from
${\mathcal N}_d^{\ast}$ (here we omit the details).

$\bullet$ Coarse subspace

With the vertex set $\mathcal{N}_d$ defined above, we can construct a coarse
subspace as in  \cite{Cai1995} and \cite{ChanZ1994}. Let ${\mathcal T}_d$ denote the auxiliary coarse mesh
generated by the vertex set $\mathcal{N}_d$, and let $\hat{V}_d(\Omega)$ denote the finite element space induced by ${\mathcal T}_d$.
For the case with irregular subdomains, the space $\hat{V}_d(\Omega)$ is not a subspace of $V_h(\Omega)$ yet, except for some very particular examples. Let $\Pi_h : \hat{V}_d(\Omega) \rightarrow V_h(\Omega)$
be the standard interpolation operator. Then the coarse space $V_d(\Omega) \subset V_h(\Omega)$ can be defined by:
\begin{equation}
V_d(\Omega)=\Big\{\Pi_h  \v: ~ \v \in \hat{V}_d(\Omega) \Big\}. \label{coarse_aux}
\end{equation}
The numerical results in section 5 will show that this coarse space is practical for elliptic-
type problems in the irregular partition situation.

$\bullet$ Vertex-related subdomains

For each vertex $\vv\in \mathcal{N}_d$, we consider all the subdomains that contain $\vv$ as their common vertex. Let $\Xi_{\vv}$ denote the set of all the vertices of
these subdomains except $\vv$ itself.
We use the software in Tetgen \cite{Tetgen} to generate a polyhedron $\hat{\Omega}_{\vv}$, which is the convex hull of the vertexes in $\Xi_{\vv}$. The vertex $\vv$ can be
regarded as the ``center" of the polyhedron $\hat{\Omega}_{\vv}$. Then we shrink this polyhedron to new one $\hat{\Omega}^{half}_{\vv}$ with half size of $\hat{\Omega}_{\vv}$, but keep the ``center" $\vv$.
Now we define the $\vv$-related subdomain $\Omega^{half}_{\vv}$ as the union of all the fine elements whose vertexes are contained in the contracted polyhedron $\hat{\Omega}^{half}_{\vv}$.
Here some auxiliary vertex subdomains $\Omega^{half}_{\vv}$ perhaps need to be extended with several fine element layers such that the union of all the
local interfaces $\Omega^{half}_{\vv}\cap\Gamma$ form an open covering of $\Gamma$.

\subsection{Comparisons between the proposed method and some existing methods}
In this subsection, we further investigate the proposed preconditioner $B$ and compare it with some existing preconditioners

$\bullet$ On the coarse solver

According to the explanations in Introduction, the coarse solver $A_d$ described in Subsection 3.2 is the simplest and almost the cheapest one in the existing non-overlapping DDMs
for three-dimensional problems (see Introduction for the details). More importantly, the construction of $A_d$ is unified and does not depend on the considered models.

$\bullet$ On the interface solvers

When the proposed coarse solver $A_d$ is used, cheap ``edge" solvers
(and ``face" solvers) need to be designed. Comparing \cite{HuSZ2011} with \cite{HuS2010} (see also \cite{DryjaSW1994} and \cite{XuZ1998}), we know that the constructions of the existing ``edge" solvers are based on
estimates of the norms induced from the interface operators restricted on the edges and so depend on the considered models. In the proposed preconditioner $B$, the construction of the vertex-related solvers $B_{\vv}$
(which play the role of interface solver) is unified and independent of the bilinear $\A(\cdot,\cdot)$.

Constructing interface solver by overlapping domain decomposition was also considered in the vertex space method (see \cite{Smith1992}) and the interface overlapping additive Schwarz method (see Subsection 7.2 of \cite{XuZ1998}). But the spaces $V^{\bot}_h(\o_{\vv}^{half})$ introduced in Subsection 3.2 have essential difference from the local interface spaces proposed in \cite{Smith1992}) and \cite{XuZ1998}, where exact harmonic extensions in all $\Omega_k$ were required (moreover, large overlap was emphasized in \cite{XuZ1998}). Thanks to this difference, we can define the cheaper (inexact) local interface solvers $B_{\vv}$.

$\bullet$ Comparisons with the BDDC method

Undoubtedly, the BDDC method is one of the most interesting non-overlapping domain decomposition methods. The BDDC method was first introduced in \cite{Dohrmann2003},
and then was studied and developed in many papers, see, for example, \cite{DohrmannW2015}-\cite{Dry2007}, \cite{KimT2009}-\cite{LiW2007}, \cite{Mandel2003}-\cite{Mandel2005} and \cite{VeigaPSWZ2014}.
The main idea of the BDDC method is to build a coarse component by the weighted sum of functions that minimize discrete local energies subject
to certain primal constraints across the subdomain interface. How to choose suitable constraints, which heavily depends on the considered models, is the key technique in the BDDC method, as in FETI-DP method.
The BDDC algorithms solve linear systems of primal unknowns, in contrast to FETI-DP algorithms.

Although the proposed method and the BDDC method were designed based on different ideas, they possess the following common merits: (1) the coarse matrix has a favorable sparsity pattern since two coarse dofs
are coupled in the coarse stiffness matrix only if both dofs appear together in at least one substructure;
(2) all the subproblems to be solved preserve the positive definiteness of the original problem; (3) the coarse stiffness matrix can be preconditioned more easily since it has almost the same structure as the original
stiffness matrix. The proposed method and the BDDC method have two respective advantages: for the proposed method, the designs of the coarse solver and local solvers are unified and do not depend on the considered model, 
and the basis functions in the coarse space can be directly obtained, without solving minimization problems; for
the BDDC method, no coarse mesh is involved so it is particularly attractive if any coarse mesh cannot directly generate a subspace of the solution, also the BDDC preconditioner in general has very fast convergence.



\section{Analysis on the convergence of the preconditioner}

\subsection{A general result}

It is known that, when the simplest coarse solver $A_d$ is used to the construction of substucturing preconditioners, the condition number of the resulting preconditioned system is not nearly optimal except for some particular cases (the coefficients in the considered equation have no large jump or the subdomains have no internal cross-point). Fortunately, this unsatisfactory condition number has no large influence on the convergence rate of the PCG iteration for solving the underlying system, provided that the number of the ``bad" eigenvalues of the preconditioned system is small (refer to \cite{XuZ2008}).

Let $\tilde{V}_h(\Omega)$ be a subspace of $V_h(\Omega)$, and let $m_0=dim(V_h(\Omega))-dim(\tilde{V}_h(\Omega))$. Assume that the number $m_0$ is small and is independent of the mesh size $h$
and the subdomain size $d$. We use $\lambda_{m_0+1}(B^{-1}A)$ to denote the minimal eigenvalue of the restriction of $B^{-1}A$ on the subspace $\tilde{V}_h(\O)$, and define $\kappa_{m_0+1}(B^{-1}A)$
as the {\it reduced condition number} of $B^{-1}A$ associated with the subspace $\tilde{V}_h(\O)$. Namely,
$$\kappa_{m_0+1}(B^{-1}A)={\lambda_{\max}(B^{-1}A)\over \lambda_{m_0+1}
(B^{-1}A)}. $$
Then the convergence rate of the PCG method with the preconditioner $B$ for solving the system (\ref{eq:2.3}) is determined by the {\it reduced condition number} $\kappa_{m_0+1}(B^{-1}A)$ (see \cite{XuZ2008} for the details).
In this subsection, we give a general result for the estimate of $\kappa_{m_0+1}(B^{-1}A)$.

For ease of notation, the symbols $\stl,~\stg$
and $\ste$ will be used in the rest of this paper. That $x_1\stl
y_1,~x_2\stg y_2$ and $x_3\ste y_3$, mean that $x_1\leq
C_1y_1,~x_2\geq c_2y_2$ and $c_3x_3\leq y_3\leq C_3x_3$ for some
constants $C_1,~c_2,~c_3$ and $C_3$ that are independent of $h$ and
$d$.

Let $I_h^{\vv}: V_h(\Omega)\rightarrow V^0_h(\Omega^{half}_{\vv})$ denote a local ``interpolation-type" operator such that, for any $\w_h\in V_h(\Omega)$, the unit decomposition condition holds, i.e., $\sum_{\vv}I_h^{\vv}\w_h\equiv \w_h$ on $\Omega$.

We define the ``extension-type" operator $E^{\bot}_{\vv}: V^0_h(\Gamma^{half}_{\vv})\rightarrow V_h^{\bot}(\Omega_{\vv}^{half})$ such that, for any $\w_h\in V^0_h(\Omega_{\vv})$, the function $E^{\bot}_{\vv}\w_h$
possesses the same degrees of freedom as $\w_h$ in $\Gamma^{half}_{\vv}$. By the definition of $V_h^{\bot}(\Omega_{\vv}^{half})$, the function $E^{\bot}_{\vv}\w_h$ vanishes in the exterior of $\Omega_{\vv}^{half}$
and is discrete ${\mathcal A}$-harmonic in $\Omega_{\vv,k}^{half}$ for each $k\in\Lambda_{\vv}$.

In the rest of this paper, let $||\cdot||_A $ denote the norm induced by the inner-product $\mathcal{A}(\cdot,\cdot)$.

\begin{theorem} Assume that there exists an operator $\Pi_d: V_h(\Omega)\rightarrow V_d(\Omega)$ such that
\begin{equation}
\|\Pi_d\v_h\|^2_A\leq C(h,d)\|\v_h\|^2_A \label{4.new-new0}
\end{equation}
and
\begin{equation}
\sum\limits_{\vv}\|E^{\bot}_{\vv}(I_h^{\vv}(\v_h-\Pi_d\v_h))\|^2_A\leq C(h,d)\|\v_h\|^2_A \label{4.new-new1}
\end{equation}
for any $\v_h\in \tilde{V}_h(\Omega)$. Here the positive number $C(h,d)$ may weakly depends on $h$ and $d$. Then we have $\kappa_{m_0+1}(B^{-1}A)\stl C(h,d)$.
\end{theorem}
\no{\it Proof}. Notice that the bilinear form ${\mathcal A}(\cdot,\cdot)$ possesses the local property, i.e., when $\v,\w\in V_h(\Omega)$ have disjoint support sets, we have ${\mathcal A}(\v,\w)=0$.
Then we can prove $\lambda_{\max}(B^{-1}A)\stl 1$ in the standard manner. It suffices to estimate $\lambda_{m_0+1}(B^{-1}A)$.

For any $\v_h\in \tilde{V}_h(\Omega)$, define $\v_d \in V_d(\Omega)$ as $\v_d =\Pi_d\v_h$ and set $\tilde{\v}_h = \v_h-\v_d$.
With this $\tilde{\v}_h$, we define the function $\tilde{\v}_{\vv}\in V_h^0(\Omega_{\vv})$ such that $\tilde{\v}_{\vv}=I_h^{\vv}\tilde{\v}_h$ on $\Gamma^{half}_{\vv}$ and $\tilde{\v}_{\vv}$
is discrete ${\mathcal A}$-harmonic in each subdomain $\Omega_k$. Then $\tilde{\v}_{\vv}\in V_h^{\bot}(\Omega_{\vv})$. Define
\begin{equation}
 \tilde{\v}_{k}^0 = (\tilde{\v}_h -
 \sum\limits_{\vv \in \mathcal{N}_d}\tilde{\v}_{\vv} )|_{\Omega_k}.\label{4.new-new2}
\end{equation}
Since the operators $I_h^{\vv}$ satisfy the unit decomposition condition, we have $ \tilde{\v}_{k}^0 \in V_h^0(\Omega_k)$.
It follows by (\ref{4.new-new2}) that
\begin{equation}
 \sum\limits_{k = 1}^{N}\tilde{\v}_{k}^0 = \tilde{\v}_h -
 \sum\limits_{\vv \in \mathcal{N}_d}\tilde{\v}_{\vv}.
 \label{eq:interdefine}
\end{equation}
Then we build the decomposition
\begin{equation}
\v_h = \v_d + \sum\limits_{k = 1}^{N}\tilde{\v}^0_{k} + \sum\limits_{\vv \in \mathcal{N}_d}
\tilde{\v}_{\vv}. \label{4.new-new5}
\end{equation}
We need only to verify the stability of the decomposition.

By the definition of $B_{\vv}$, we have
$$ (B_{\vv}\tilde{\v}_{\vv},\tilde{\v}_{\vv})=\|E^{\bot}_{\vv}(I_h^{\vv}\tilde{\v}_h)\|^2_A. $$
This, together with the assumption (\ref{4.new-new1}), leads to
\begin{equation}
 \sum\limits_{\vv}(B_{\vv}\tilde{\v}_{\vv},\tilde{\v}_{\vv})\leq C(h,d)||\v_h||^2_A.\label{4.new-new4}
\end{equation}
Finally, using (\ref{4.new-new5}), together with (\ref{4.new-new0}) and (\ref{4.new-new4}), yields
\bn
&&(A_d\v_d,\v_d) + \sum\limits_{k=1}^N(A_k \tilde{\v}^0_{k}, \tilde{\v}^0_{k})+
 \sum\limits_{\vv \in \mathcal{N}_d}(B_{\vv}\tilde{\v}_{\vv},\tilde{\v}_{\vv}) \cr
&& \stl C(h,d)||\v_h||^2_A, \quad\forall \v_h \in \tilde{V}_h(\Omega),
 \label{eq:proof}
 \en
which implies that
$$ \lambda_{m_0+1}(B^{-1}A)\geq 1/C(h,d). $$
Now the desired result can be obtained directly.
\hfill$\Box$

\subsection{Result for linear elasticity problems}

In this subsection, we try to estimate the positive $C(h,d)$ in Theorem 4.1 for linear elasticity problems.
As to Maxwell's equations, the proof is very complicated, which beyond the goal of this article.

Let's consider the linear elasticity problem:
\be \begin{cases}
-\sum\limits_{j=1}^3\frac{\partial \sigma_{ij}}{\partial x_{j}}(\bm{u})=f_{i},\q in~~\O \\
 \q\q\q\q\q \bm{u}=0, \q on~~\partial\Omega \end{cases}\l{eq:5.1}
\ee
where $\bm{f}=(f_1~f_2~f_3)^{T}$ is an internal volume force, e.g. gravity (cf. \cite{ChenH2011}). The linearized strain tensor is defined by
\[\varepsilon = \varepsilon(\bm{u}) = [\varepsilon_{ij} = \frac{1}{2}(\frac{\partial
u_{i}}{\partial x_{i}} +\frac{\partial u_{j}}{\partial x_{i}})], \]
and
\[ \sigma_{ij}(\bm{u}) := \lambda\delta_{ij}div\bm{u} + 2\mu\varepsilon_{ij}. \]
where $\lambda$ and $\mu$ are the $Lam\acute{e}$ parameters (cf. \cite{re25}), which are positive
functions.

The subspace $H_{0}^1(\Omega)\subset H^1(\Omega)$ is the set of
functions having the zero trace on $\partial \Omega$. We introduce the vector
value Sobolev space $(H_{0}^1(\Omega))^3$.
Concerning the variational problem (\ref{eq:2.1}), we have $V(\Omega):=[H_{0}^1(\Omega)]^3$ and

\[ \A(\bm{u},\bm{v}) = \int_{\Omega}(2\mu\varepsilon(\bm{u}):\varepsilon(\bm{v}) + \lambda div\bm{u}\cdot div\bm{v})d\x, \]
\[ \langle \bm{F},\bm{v} \rangle = \int_{\Omega}\bm{f}\cdot\bm{v}d\x \]
where
\[\varepsilon(\bm{u}):\varepsilon(\bm{v}) := \sum_{i,j=1}^n \varepsilon_{ij}(\bm{u})\varepsilon_{ij}(\bm{v}).\]

Let $R(K)$ be a subset of all linear polynomials on the element
$K$ of the form:
$$ R(K)=\Big\{\bm{A}\cdot\x+\bm{C};~ \bm{A}\in \mathds{R}^{3\times3}, \bm{C} \in \mathds{R}^3, ~\x\in K\Big\} . $$

To our knowledge, there seems no work to analyze nearly optimal substructuring preconditioners for three-dimensional problems with irregular subdomains, which bring particular difficulties. 
Thus, here we only consider the case with regular subdomains. Assume that $\O$ can be written as the union of polyhedral subdomains $D_1$, $\cdots$, $D_{N_0}$, 
such that $\lambda(\x) =\lambda_r$ and $\mu(\x) = \mu_r$
for $x\in D_r$, with every $\lambda_r$ and $\mu_r$ being a positive constant. In applications, $N_0$ is a {\it fixed} positive integer, so the diameter of each $D_r$ is $O(1)$.
For the analysis, we assume that every subdomain $\Omega_k$ is a polyhedron.
It is certain that the subdomains $\Omega_k$ should satisfy the condition: each $D_r$ is the union of some subdomains in $\{\Omega_k\}$.

The null space $ker(\varepsilon)$ is the space of rigid body motions. In three
dimensions, the corresponding space is spanned by three translations
$$
\r_1:=
\begin{bmatrix}
1\\
0\\
0
\end{bmatrix}
,\quad \r_2:=
\begin{bmatrix}
0 \\
1 \\
0
\end{bmatrix},
\quad\r_3:=
\begin{bmatrix}
0 \\
0 \\
1
\end{bmatrix},
$$
and three rotation
$$
\r_4:=
\begin{bmatrix}
0\\
x_3\\
-x_2
\end{bmatrix}
,\quad \r_5:=
\begin{bmatrix}
-x_3 \\
0 \\
x_1
\end{bmatrix},
\quad\r_6:=
\begin{bmatrix}
-x_2 \\
x_1 \\
0
\end{bmatrix}.
$$

Let $\Lambda=\{k:~\partial
D_k\cap\partial\O=\emptyset\}$ denote the index set of the
subdomains $D_1,\cdots,D_{N_0}$ that do not touch the boundary of
$\O$, and set
$$ \tilde{V}_h(\O)=\{\v_h\in V_h(\O):~\sum\limits_{i=1}^6 |\int_{D_k}\r_i \cdot \v_h
d\x|=0,~k\in \Lambda\}.$$
Let $n_0$ denote the number of the subdomains that do not touch $\partial\Omega$, i.e., $n_0=dim(\Lambda)$, and set $m_0 = 6n_0$.

\begin{theorem} For the linear elasticity problems described above, we have
\begin{equation}
\kappa_{m_0+1}(B^{-1}A)\stl
\log(1/d)\log^2(d/h). \l{eq:theorem}
\end{equation}
When the coefficients $\mu(x)$ and $\lambda(x)$ have no large jump across the interface $\Gamma$, or there is no cross-point in the distribution of the jumps of the coefficients, the factor $\log(1/d)$ in the above inequalities
can be removed.
\end{theorem}

In order to prove this result, we need several auxiliary results.
In the following we use $D\subset \O$ to denote a generic polyhedron in $D_1,\cdots,D_{N_0}$.
It is clear that
\begin{equation}
||\varepsilon(\v)||_{0,D} \leq ||\nabla \v||_{0,D}\quad\mbox{and}\quad
||div(\v)||_{0, D} \leq ||\nabla \v||_{0,D}, \forall \v \in [H^1(D)]^3.
\label{eq:new}
\end{equation}

The following lemma can be obtained by Korn inequality and the result in \cite{KlawWD2000} and \cite{Necas}.
\begin{lemma}
There exist positive constants  $c_0$ and $C_0$, such that
$$ c_0||\nabla \v||_{0,D} \leq ||\varepsilon (\v)||_{0,D} \leq C_0||\nabla \v||_{0,D} $$
for any function $\v\in [H^1 (D)]^3$, which satisfies either $(\v, \r)_{L_2(D)} = 0$ for each $\r\in ker(\varepsilon)$ or $\v$ vanishes on a face of $D$.
 \end{lemma}

\hfill$\Box$

 Define the weighted norm
   $$||\v||_{1,\O_k} = (|\v|^2_{1,\O_k} + d^{-2}||\v||^2_{0,\Omega_k})^{1/2}$$

We assume that there exists constant $c_0, C_0$, such that
$$c_0 \lambda_r \leq \mu_r \leq C_0 \lambda_r,
 \quad \forall D_r\subset \O, \quad  r = 1, \cdots, D_{N_0}.$$
Define the weighted $L^2$-inner product and the weighted $H^1$-inner product as follows:
   $$(\u, \v)_{L^2_{\lambda}(\Omega)} = \sum\limits_{r = 1}^{N_0}\lambda_r\int_{D_r} \u \cdot \v d\x,
   \quad\u,   \v \in [L^2(\Omega)]^3.
   $$
   $$(\u, \v)_{H^1_{\lambda}(\Omega)} = \sum\limits_{r = 1}^{N_0}\lambda_r\int_{D_r} \nabla\u \cdot \nabla\v d\x,
   \quad\u,   \v \in [H^1_0(\Omega)]^3.
   $$
Let $||\cdot||_{L^2_{\lambda}(\Omega)}$ and $|\cdot|_{H^1_{\lambda}(\Omega)}$
denote, respectively, the norm and the semi-norm induced by the inner product
$(\cdot,\cdot)_{L^2_{\lambda}(\Omega)}$ and $(\cdot,\cdot)_{H^1_{\lambda}(\Omega)}$. For convenience, define
$$||\v||_{H^1_{\lambda}(\Omega)} = (|\v|^2_{H^1_{\lambda}(\Omega)} + d^{-2}||\v||^2_{L^2_{\lambda}(\Omega)})^{1/2}$$

Let $Q_d^{\lambda}: [L^2(\O)]^3\ri V_d(\O)$ be the weighted $L^2$
projections defined by
\begin{equation}
  \label{wproject}
   (Q_d^{\lambda} \v, \w)_{L^2_\lambda(\Omega)}=(\v, \w)_{L^2_\lambda(\Omega)}, \quad
            \forall \v\in [L^2(\Omega)]^3, \w\in V_d(\O).
\end{equation}

It is clear that
$$\tilde{V}_h(\O) \subset \{\v_h\in V_h(\O):~\int_{D_k}\v_h
d\x={\bf 0},~k\in \Lambda\}.$$
Then, by the result given in \cite{XuZ2008}, we have
\begin{lemma}
The weighted $L^2$ projection $Q_d^{\lambda}$ satisfies
\begin{equation}
\|(Q_d^{\lambda}-I)\v\|^2_{L^2_{\lambda}(\Omega)} \stl d^2\log(1/d)|\v|^2_{H^1_{\lambda}(\Omega)}, \q \forall \v \in
\tilde{V}_h(\O) \l{project1}
\end{equation}
and \be |Q_d^{\lambda}\v|^2_{H^1_{\lambda}(\Omega)}\stl\log(1/d)|\v|^2_{H^1_{\lambda}(\Omega)},\q\forall \v\in
\tilde{V}_h(\O).\l{project2} \ee
\end{lemma}
\hfill$\Box$

\begin{remark}
When the coefficient $\lambda(\x)$ has no jump across the interface
$\G$, or there is no cross-point in the distribution of the jumps of
the coefficient, the factor $\log(1/d)$ in the inequalities
(\ref{project1}) and (\ref{project2}) can be removed.
\end{remark}

 The following lemma is a direct consequence of Lemma 4.1 and (\ref{eq:new}).

\begin{lemma} The following norm equivalence holds
\begin{equation}
  ||\v||^2_{A} \ste |\v|^2_{H^1_{\lambda}(\Omega)}, \quad \forall \v \in \tilde{V}_h(\Omega).
  \label{eq:equalnorm}
  \end{equation}
\end{lemma}

\hfill$\Box$

Next we present a stability result of discrete harmonic functions in some $\Omega_{\vv,k}^{half}$.
In order to simplify the analysis, we assume that all our subdomains $\O_k$ are cubes and we only consider a specific way for the construction of
$\Omega_{\vv}^{half}$.

For each $\vv \in \mathcal{N}_d$, we choose an auxiliary cube $G_{\vv}$ with the
size $d$ and the center $\vv$.
Define $\Omega_{\vv}^{half}$ as the union of all the elements, each of
which has at least one vertex located in the cube or just touching
the boundary of this cube. Notice that $G_{\vv}$ may be sightly smaller than $\Omega_{\vv}^{half}$ (we consider only the case with small overlap).
In this case, the size of
$\O_{\vv,k}^{half}$ is $d/2$ for each $k\in\Lambda_{\vv}$. Notice that $\O_{\vv,k}^{half}$ is not a polyhedron, i.e., it has a quite
irregular boundary, so the desired stability result can not be proved in the standard way. It is easy to see that $\O_{\vv,k}^{half}$
is a Jones domain, whose definition and properties can be found in \cite{Jones1981} and \cite{Klawonn2008}. Intuitively, a Jones domain means that
it satisfies twisted cone condition  and can't be too oblate relative to its diameter.

Let $\mathcal{E}_{ij}$ be the common face between two neighboring subdomains $O_i$ and $O_j$. Define
$$ V_h^{ij}(O_k)= \{ v_h\in V_h(O_k): v_h(\x) = 0~ at~ all~ nodes~ of~ \partial O_k\backslash\mathcal{E}_{ij}\}\quad(k=i,j).$$
Before giving the stability result, we recall an ``interface" extension lemma proved in \cite{Klawonn2008}.
\begin{lemma} \cite{Klawonn2008} Assume that $O_i$ is a domain with a complement which is a Jones domain.
 Then, there exists an extension operator
            $$E^h_{ji}:V_h^{ij}(O_j) \rightarrow V_h^{ij}(O_i)$$
such that
$$(E^h_{ji}v_h)|_{\Gamma_{ij}}=v_h|_{\Gamma_{ij}}\quad\mbox{and}\quad |E^h_{ji}v_h|_{H^1(O_i)} \stl|v_h|_{H^1(O_j)},\quad\forall v_h\in V_h^{ij}(O_j).$$
\end{lemma}

\hfill $\Box$


To distinguish with
the elastic harmonic extension,  we define the discrete harmonic space associated with Laplace-type operator
$$ V^{H}_h(\Omega)=\{\v\in V_h(\o):~~(\nabla\v,\nabla\w)=0,~\forall\w\in
V_h^{0}(\o_{k}),~~k= 1 \cdots N\}.$$
For each $\vv\in {\mathcal N}_d$, define the vertex-related local harmonic spaces
$$ V^{H}_h(\Omega_{\vv})=\{\v\in
V^{H}_h(\o):~~\mbox{the~trace~of}~\v~\mbox{belongs~to}~V_h^0(\Gamma_{\vv}^{half})\}\subset
V_h^0(\Omega_{\vv})$$
and
$$ V^{H}_h(\o_{\vv}^{half})=\{\v\in V^0_h(\o_{\vv}^{half}):~~(\nabla\v,
\nabla\w)=0,~\forall\w\in V_h^{0}
(\o_{\vv,k}^{half})~~\mbox{with}~~k\in\Lambda_{\vv}\}.$$

\begin{lemma} Let $\u^{half}_{\vv} \in V_h^{H}(\Omega_{\vv}^{half})$ and $\u_{\vv}
\in V_h^{H}(\Omega_{\vv})$. Assume that the two functions satisfy $\u_{\vv} = \u^{half}_{\vv}~ on~
\Gamma_{\vv}^{half}$. Then we have
\begin{equation}
 |\u^{half}_{\vv}|_{H^1_{\lambda}(\Omega)} \stl|\u_{\vv}|_{H^1_{\lambda}
(\Omega)}.\label{eq:crossequal}
\end{equation}
\end{lemma}
\no{\it Proof}. For $k \in
\Lambda_{\vv}$,  we define $\u^{half}_{\vv,k} = \u^{half}_{\vv}|_{\Omega_k}$ and
$\u_{\vv,k} = \u_{\vv}|_{\Omega_k}$. It suffices to prove that
\begin{equation}
 |\u^{half}_{\vv,k}|_{1,\Omega_k} \stl |\u_{\vv,k}|_{1,\Omega_k},\quad\forall k\in \Lambda_{\vv}.
 \label{eq:crossequal1}
 \end{equation}
By the triangle inequality, we have
\begin{equation}
|\u^{half}_{\vv,k}|_{1,\Omega_k} \leq |\u_{\vv,k}|_{1,\Omega_k}
+|\u^{half}_{\vv,k} - \u_{\vv,k}|_{1,\Omega_k}.
\label{eq:crossequal2}
\end{equation}
Set $\Omega^{\partial}_{\vv,k} = \Omega_k\backslash\Omega_{\vv,k}^{half}$ and let $\ff_{\vv,k}$ denote the common interface of $\Omega^{\partial}_{\vv,k}$
and $\Omega_{\vv,k}^{half}$. Since $\u^{half}_{\vv,k}$ vanishes on $\ff_{\vv,k}$, it can be naturally
extended into $\Omega^{\partial}_{\vv,k}$ by zero. Let $\tilde{\u}^{half}_{\vv,k}\in V_h(\Omega_k)$ denote the resulting extension, and define $\u^{\partial}_{\vv,k} =\u_{\vv,k} - \tilde{\u}^{half}_{\vv,k}$.
Then $\u^{\partial}_{\vv,k} |_{\Omega^{\partial}_{\vv,k}} = \u_{\vv,k}$ and $\u^{\partial}_{\vv,k}$ vanishes on $\partial\Omega^{\partial}_{\vv,k}
\backslash\ff_{k,\vv}$.
Using Lemma 4.4, there exists an extension $\tilde{\u}^{\partial}_{\vv,k}$ of $\u^{\partial}_{\vv,k}$ such that $\tilde{\u}^{\partial}_{\vv,k}\in V_h(\Omega_{\vv,k}^{half})$ and
$\tilde{\u}^{\partial}_{\vv,k}=\u^{\partial}_{\vv,k}$ on $\ff_{\vv,k}$. Moreover, the extension $\tilde{\u}^{\partial}_{\vv,k}$ satisfies
\begin{equation}
|\tilde{\u}^{\partial}_{\vv,k}|_{1,\Omega_{\vv,k}^{half}} \stl |\u^{\partial}_{\vv,k} |_{1,\Omega^{\partial}_{\vv,k}}\quad\mbox{and}\quad\tilde{\u}^{\partial}_{\vv,k} = \u^{\partial}_{\vv,k} ~on ~
 \partial\Omega_{\vv,k}^{half}.
\label{eq:crossequal3}
\end{equation}

Notice that $\u^{\partial}_{\vv,k}$ is Laplace-type discrete harmonic in the interior of $\Omega_{\vv,k}^{half}$, by (\ref{eq:crossequal3}) we obtain
\begin{equation}
 |\u_{\vv,k}^{\partial}|_{1,\Omega_{\vv,k}^{half}} \leq |\tilde{\u}^{\partial}_{\vv,k}
 |_{1,\Omega_{\vv,k}^{half}} \stl|\u^{\partial}_{\vv,k} |_{1,\Omega^{\partial}_{\vv,k}}
 = |\u_{\vv,k} |_{1,\Omega^{\partial}_{\vv,k}}.
 \label{eq:crossequal4}
\end{equation}
Combing this with (\ref{eq:crossequal2}) leads to  (\ref{eq:crossequal1}).
Then (\ref{eq:crossequal}) follows by (\ref{eq:crossequal1}).
\hfill $\Box$

\begin{lemma}
Let $\u^{\perp} \in V_h^{\perp}(\Omega_{\vv})$,  $\u^H \in V_h^{H}(\Omega_{\vv})$ and
$\u^{\perp} = \u^H ~on~ \Gamma_{\vv}^{half}$. Then we have
\begin{equation}
(B_{\vv} \u^{\perp}, \u^{\perp}) \stl|\u^H|_{H^1_{\lambda}}
\label{eq:bvnormchange}
\end{equation}
\end{lemma}
\no{\it Proof}. Define $\u^{half} \in V^{\perp}(\O_{\vv}^{half})$ such that $\u^{half} = \u^{\perp} ~on~ \Gamma_{\vv}^{half}$. According to the definition
of $B_{\vv}$, we know that
\begin{equation}
(B_{\vv} \u^{\perp}, \u^{\perp}) = ||\u^{half}||_A^2.
\label{eq:bvnormchange1}
\end{equation}
Let $\tilde{\u}^{half} \in V^{H}(\O_{\vv}^{half})$ and satisfy $\tilde{\u}^{half} = \u^{half} ~on~ \Gamma_{\vv}^{half}$. For each
$k \in \Lambda_{\vv}$,~the function $\u^{half}$ is elastic-type discrete harmonic and $\tilde{\u}^{half}$ is Laplace-type discrete harmonic in the interior of $\Omega^{half}_{\vv,k}$. Then, by (\ref{eq:new}), we have
\begin{equation}
||\u^{half}||^2_{A} \stl |\tilde{\u}^{half}|
^2_{H^1_{\lambda}
(\Omega)}.
\label{eq:bvnormchange2}
\end{equation}
Using Lemma 4.5, we know that
\begin{equation}
 |\tilde{\u}^{half}|^2_{H^1_{\lambda}(\Omega)} \stl|\u^{H}|^2_{H^1_{\lambda}
(\Omega)}.
\label{eq:bvnormchange3}
\end{equation}
Then (\ref{eq:bvnormchange}) follows by (\ref{eq:bvnormchange1}),
(\ref{eq:bvnormchange2}) and (\ref{eq:bvnormchange3}).
\hfill $\Box{}$

For a vector-valued function $\v =(v_1~ v_2~ v_3)^T \in (H^1(\Omega_k))^3$, we define the $H^{1/2}$-norm
$$  ||\v||_{1/2,\partial\Omega_k}=( |\v|^2_{1/2,\partial\Omega_k}
   + d^{-1}||\v||^2_{0,\partial\Omega_k})^{1/2}~~\mbox{with}~~|\v|^2_{1/2,\partial\Omega_k}= \sum\limits_{i = 1}^3
|v_i|^2_{1/2,\partial\Omega_k}.$$
For a face $F$ of $\partial\O_k$ and $\v_h\in V^0_h(F)$, let $\tilde{\v}_h\in V_h(\partial\O_k)$ denotes the zero extension of $\v_h$. Define the norm
$$
\|\v_h\||^2_{H^{1/2}_{00}(F)}=|\tilde{\v}_h|^2_{1/2,~\partial\O_k}.$$

For a given subset $K\subset\O$, define a restriction operator $I_K^0:V_h(\O)\rightarrow V_h^0(K)$ as follows: $(I^0_K\v)(x) =\v(x)$ for any $x\in K\cap
\mathcal{N}_h$
and $(I^0_K\v)(x)={\bf 0}$ for $x\in\Omega\backslash K$. Similarly, we can define $I^0_K:V_h(\Gamma)\rightarrow V_h^0(K)$ for a subset $K$ of the interface $\Gamma$.

Let $\vv$ be a given vertex. For a face $\ff$ that has $\vv$ as one of its vertices, let $\ff^{in}_{\vv}$ denote the intersection of $\ff$ and $\Gamma_{\vv}^{half}$.
In the following we prove an extension result on the norm $H^{{1/ 2}}_{00}(\ff^{in}_{\vv})$.

\begin{figure}[ht]
\centering
\includegraphics[width=9cm]{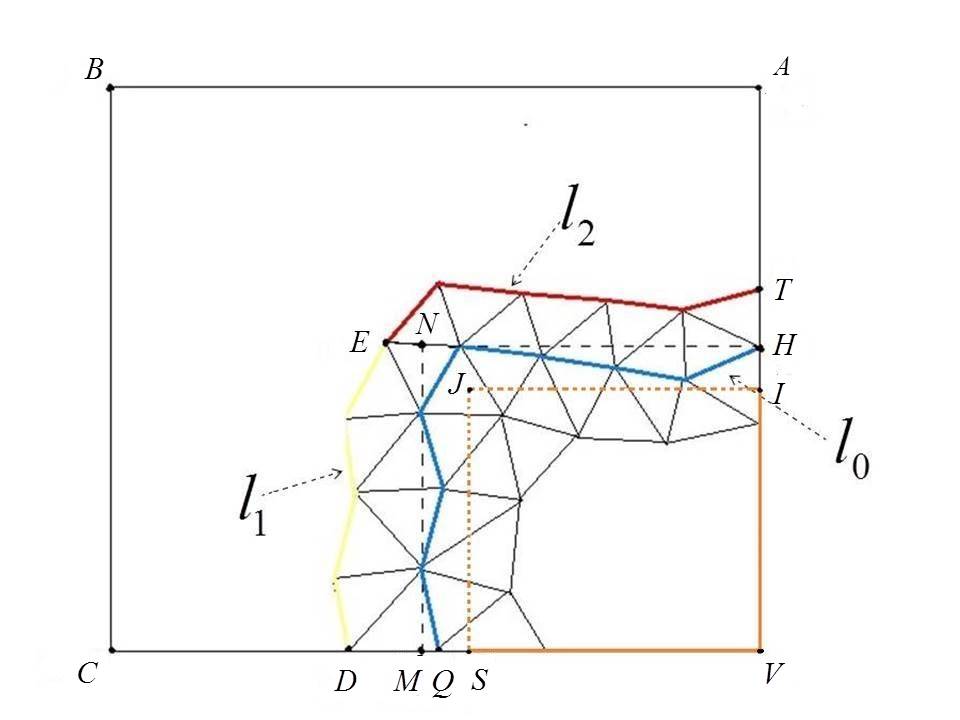}
\caption{the relation between $\ff$ and $\ff^{in}_{\vv}$}
\label{fig:triangle}
\end{figure}

For convenience, we show $\ff$ and $\ff^{in}_{\vv}$ in Fig. \ref{fig:triangle}, where the big square $ABCV$ denotes $\ff$. Let $l_0$ , $l_1$, $l_2$ denote the broken segments $QH$ (the blue curve), $DE$ (the yellow curve) and $ET$ (the red curve), respectively. Then $\ff^{in}_{\vv}$ is just the area surrounded by $l_1$, $l_2$, the
straight segments $\overline{DV}$ and $\overline{TV}$. In Fig. \ref{fig:triangle}, the smaller square $HNMV$ denotes the intersection of $\ff$ with the auxiliary cube $G_{\vv}$ described behind Lemma 4.3. For the proof, we define an auxiliary square $SVIJ$ (which is denoted by $\ff^{aux}_{\vv}$) such that the distance between $\overline{MN}$ and $\overline{SJ}$
is approximately equal to  $h/2$.

\begin{lemma} For a subdomain $\Omega_k$, let $\ff^{in}_{\vv}\subset\partial\Omega_k$ be defined above. Then we have the face inequality
\begin{equation}
||I^0_{\ff^{in}_{\vv}} \v_h||_{H^{1/2}_{00}(\ff^{in}_{\vv})} \stl
\log(d/h)||\v_h||_{1,\O_k},\quad\forall \v_h \in \ V_h(\O_k).
\label{eq:face0}
\end{equation}
\end{lemma}
\no{\it Proof}. As in Lemma 3.1 of \cite{DohrmannW2015}, we can prove the ``edge" inequality
\begin{equation}
||\v_h||_{0,l_{i}} \stl \log^{1/2}(d/h)||\v_h||_{1,\O_k}, \quad i =0, 1, 2.
\label{eq:edges}
\end{equation}

 It is known that
 \begin{equation}
 ||I^0_{\ff^{in}_{\vv}} \v_h||^2_{H^{1/2}_{00}(\ff^{in}_{\vv})}\ste
|I^0_{\ff^{in}_{\vv}} \v_h|^2_{1/2, \ff^{in}_{\vv}}+
\int_{\ff^{in}_{\vv}}\frac{|I^0_{\ff^{in}_{\vv}} \v_h|^2}{dist(\x,\partial\ff^{in}_{\vv})}
ds(\x).
\label{eq:define}
\end{equation}
It follows by (\ref{eq:edges}) that
\bn
|I^0_{\ff^{in}_{\vv}} \v_h|_{1/2, \ff^{in}_{\vv}} &\stl&
  |\v_h|_{1/2, \ff^{in}_{\vv}} + \sum\limits_{i= 1}^{2}||\v_h||_{0,l_{i}} + ||\v_h||_{0,\overline{DV}} + ||\v_h||_{0,\overline{GV}}  \cr
  &\stl& |\v_h|_{1/2,\partial\O_k} + \log^{1/2}(d/h)
  ||\v_h||_{1,\O_k} \cr
  &\stl&  \log^{1/2}(d/h)||\v_h||_{1,\O_k}.
  \label{eq:face1}
\en
Let $\ff^{aux}_{\vv}$ denote the auxiliary square $SVIJ$ shown in \textsc{Fig.} \ref{fig:triangle}.
Then we have
\bn
\int_{\ff^{in}_{\vv}}\frac{|I^0_{\ff^{in}_{\vv}} \v_h(\x)|^2}{dist(\x,
\partial\ff^{in}_{\vv})}ds(\x) &\ste& h^2\sum\limits_{p_i \in \mathcal{N}_h\cap\ff^{in}_{\vv}}
\frac{|\v_h(p_i)|^2}{dist(p_i, \partial \ff^{in}_{\vv})} \cr
& \stl & h^2\sum\limits_{p_i \in \mathcal{N}_h\cap l_{0}}
\frac{|\v_h(p_i)|^2}{dist(p_i, \partial \ff^{in}_{\vv})} +
h^2\sum_{\substack{p_i \in \mathcal{N}_h\cap\ff^{in}_{\vv}  \\
 p_i \notin l_0}}
\frac{|\v_h(p_i)|^2}{dist(x_i, \partial \ff^{in}_{\vv})} \cr
&\stl& \|\v_h\|^2_{0,l_0}+
\int_{\ff^{aux}_{\vv}}\frac{|\v_h(\x)|^2}{dist(\x,\partial\ff^{in}_{\vv})}ds(\x).
\label{eq:face2}
\en
As in Lemma 4.10 in \cite{XuZ1998}, we can verify that
$$
\int_{\ff^{aux}_{\vv}}\frac{|\v_h(\x)|^2}{dist(\x,\partial\ff^{in}_{\vv})}ds(\x)
\stl\log^2(d/h)||\v_h||^2_{1,\O_k}.
$$
Substituting (\ref{eq:edges}) and the above inequality into (\ref{eq:face2}), yields
$$ \int_{\ff^{in}_{\vv}}\frac{|I^0_{\ff^{in}_{\vv}} \v_h(\x)|^2}{dist(\x,
\partial\ff^{in}_{\vv})}ds(\x)\stl\log^2(d/h)||\v_h||^2_{1,\O_k}. $$
Plugging this and (\ref{eq:face1}) in (\ref{eq:define}), leads to (\ref{eq:face0}).

\hfill$\Box$
\vskip 0.2in
\no{\it Proof of Theorem 4.2}. It suffices to verify the assumptions given in Theorem 4.1. To this end,
we choose the weighted $L^2$ projector $Q_d^{\lambda}$ as the operator $\Pi_d$. Let $\tilde{\v}_h\in \tilde{V}_h(\Omega)$.

It follows by (\ref{eq:new}), (\ref{project2}) and (\ref{eq:equalnorm}) that
\begin{equation}
\|\Pi_d\v_h\|^2_A\stl\log(1/d)\|\v_h\|^2_A.\l{eq:coarsenorm}
\end{equation}
Set $\tilde{\v}_h=\v_h-\Pi_d\v_h$. Using Lemma 4.2 and Lemma 4.3, yields
\begin{equation}
||\tilde{\v}_h||^2_{H^1_{\lambda}(\O)}
\stl \log(1/d)||\v_h||^2_{A}.
\label{eq:HtoEnergy}
\end{equation}

Let $I_h^{\vv}: V_h(\Omega)\rightarrow V_h^0(\Omega_{\vv}^{half})$ be the nodal weight interpolation: for $\w_h\in V_h(\Omega)$, the function $I_h^{\vv}\w_h$ vanishes outside $\Omega_{\vv}^{half}$  and
$(I_h^{\vv}\w_h)(\x)={1\over q(\x)}\w_h(\x)$ for $\x\in \Omega_{\vv}^{half}\cap{\mathcal N}_h$, where $q(\x)$ denotes the number of the subdomains $\Omega_{\vv}^{half}$ containing $\x$. Then all the
operators $I_h^{\vv}$ satisfy the unit decomposition condition described in the front of Theorem 4.1.
Define $\tilde{\v}_{\vv}\in V_h^{\bot}(\Omega^{half}_{\vv})$ as $\tilde{\v}_{\vv}=E^{\bot}_{\vv}(I_h^{\vv}(\v_h-\Pi_d\v_h))$, where $E^{\bot}_{\vv}$ is the extension operator defined in Subsection 4.1.

In the following, we estimate $\sum_{\vv}\|\tilde{\v}_{\vv}\|^2_A$.

For the function $\tilde{\v}_{\vv} \in V_h^{\perp}(\O^{half}_{\vv})$, define $\tilde{\v}_{\vv}^{H} \in V_h^{H}(\O_{\vv})$ such that $\tilde{\v}_{\vv}^{H} = \tilde{\v}_{\vv}~ on~ \Gamma_{\vv}^{half}$.
It follows by Lemma 4.6 that
\bn
\|\tilde{\v}_{\vv}\|^2_A &\stl&
|\tilde{\v}^{H}_{\vv}|^2_{H^1_{\lambda}(\Omega)}
= \sum\limits_{k \in \Lambda_{\vv}}\lambda_k |\tilde{\v}^{H}_{\vv}|^2_{1,\Omega_k}\cr
&\stl & \sum\limits_{k \in \Lambda_{\vv}}\lambda_k |\tilde{\v}^{H}_{\vv}|^2_{1/2,\partial\Omega_k}.
\label{eq:normconvert}
\en
Thus we need only to estimate $|\tilde{\v}^{H}_{\vv}|^2_{1/2,\partial\Omega_k}$. It is clear that
\bn
 \tilde{\v}^{H}_{\vv}|_{\partial\O_k}=I^0_{\vv}\tilde{\v}^{H}_{\vv} +
\sum_{\substack{\se \in \mathcal{E}_{\vv} \\ \se \subset \Gamma_k}}I^0_{\se} \tilde{\v}^{H}_{\vv}
+\sum_{\substack{\ff \in \mathcal{F}_{\vv} \\ \ff \subset \Gamma_k}}I^0_{\ff} \tilde{\v}^{H}_{\vv}.
\label{eq:split}
\en
Notice that $I^0_{\ff}\tilde{\v}^{H}_{\vv}$ vanishes on $\ff \backslash \ff^{in}_{\vv}$. Then
\bn
&&||I^0_{\ff} \tilde{\v}^{H}_{\vv}||^2_{H^{1/2}_{00}(\ff)}\stl
||I^0_{\ff^{in}_{\vv}}(I_h^{\vv}\tilde{\v}_{h})||^2_{H^{1/2}_{00}(\ff^{in}_{\vv})} \cr
&&\stl ||I^0_{\ff^{in}_{\vv}}\tilde{\v}_{h}||^2_{H^{1/2}_{00}(\ff^{in}_{\vv})}
 +||I^0_{\ff^{in}_{\vv}}(\tilde{\v}_{h}-I_h^{\vv}\tilde{\v}_{h})||^2_{H^{1/2}_{00}(\ff^{in}_{\vv})}.
 \label{eq:facetoface}
\en
By the construction of $\Omega_{\vv}^{half}$ (see the description behind Lemma 4.3), the overlap between intersecting faces $\ff^{in}_{\vv}$ associated with two neighboring vertices contains at most two elements layer (refer to Fig. \ref{fig:triangle}). Then, by the definitions of $I_h^{\vv}$ and $I^0_{\ff^{in}_{\vv}}$, the function $I^0_{\ff^{in}_{\vv}}(\tilde{\v}_{h}-I_h^{\vv}\tilde{\v}_{h})$ vanishes at all nodes except on $l_0$ shown in Fig. \ref{fig:triangle}. In particular, when the overlap is just one element layer,
we have $I^0_{\ff^{in}_{\vv}}(\tilde{\v}_{h}-I_h^{\vv}\tilde{\v}_{h})\equiv 0$. By the inverse estimate and the discrete $L^2$ norms,
we can deduce that
$$ ||I^0_{\ff^{in}_{\vv}}(\tilde{\v}_{h}-I_h^{\vv}\tilde{\v}_{h})||^2_{H^{1/2}_{00}(\ff^{in}_{\vv})}\stl
     \|(\tilde{\v}_{h}-I_h^{\vv}\tilde{\v}_{h})||^2_{0,l_{0}}\stl\|\tilde{\v}_h||^2_{0,l_0}.$$
Substituting this into (\ref{eq:facetoface}), and using (\ref{eq:face0}) and (\ref{eq:edges}), yields
\begin{equation}
||I^0_{\ff} \tilde{\v}^{H}_{\vv}||^2_{H^{1/2}_{00}(\ff)}\stl \log^2(d/h)||\tilde{\v}_h||^2_{1,\O_k}.
 \label{eq:faceextension}
\end{equation}
In addition, using the vertex and edge lemmas in \cite{XuZ1998}, we get
$$
||I^0_{\vv}\tilde{\v}^{H}_{\vv}||_{1/2,\partial\Omega_k}
\stl ||I^0_{\vv}\tilde{\v}_h||_{1/2,\partial\Omega_k}
 \stl \log^{1/2}(d/h)||\tilde{\v}_h||_{1,\O_k}
$$
and
$$
||I^0_{\se} \tilde{\v}^{H}_{\vv}||
^2_{1/2,\partial\Omega_k}
\stl
||I^0_{\se} \tilde{\v}_h||
^2_{1/2,\partial\Omega_k}
\stl \log^{1/2}(d/h)||\tilde{\v}_h||_{1,\O_k}.
$$
By (\ref{eq:split}), together with (\ref{eq:faceextension}) and the above two inequalities, gives
$$ |\tilde{\v}^{H}_{\vv}|^2_{1/2,\partial\O_k} \stl \log^{2}(d/h)||\tilde{\v}_h||^2_{1,\O_k}. $$
Substituting this into (\ref{eq:normconvert}) and using (\ref{eq:HtoEnergy}), yields
\bq
\sum\limits_{\vv}\|\tilde{\v}_{\vv}\|^2_A&\stl&
\sum\limits_{\vv}\sum\limits_{k \in
\Lambda_{\vv}}\lambda_k \log^2(d/h)||\tilde{\v}_h||^2_{1,\Omega_k}\cr
&\stl& \log^2(d/h)||\tilde{\v}_h||^2_{H^1_{\lambda}(\O)} \cr
&\stl& \log(1/d)\log^2(d/h)\|\v_h\|^2_A.
\eq
This, together with (\ref{eq:coarsenorm}), verify the assumptions of Theorem 4.1 with $C(h,d)=\log(1/d)\log^2(d/h)$.
\hfill$\Box$

\section{Numerical Experiments}
\setcounter{equation}{0} In this section, we report some numerical results for the linear elasticity problems and Maxwell's equations to illustrate the efficiency
of the proposed substructuring methods.

In our experiments, we choose the domain $\Omega=(0,1)^3$ and define domain decomposition and finite element partition as follows.
At first, we divide the domain into $n^3$ smaller cubes $\Omega_1$, $\Omega_2\cdots\Omega_N$, which have the same length $d$ of edges, i.e., $d = 1/n$.
We require that $D\subset\o$ is just the union of some subdomains in $\{\o_k\}$, which yields the desired domain decomposition. Next, we divide each subdomain $\Omega_k$ into $m^3$
fine cubes, with the same size $h = 1/(mn)$.
Then we further divide each fine cube into 5 or 6 tetrahedrons in the standard way, then all the generated tetrahedrons constitute a partition $\mathcal{T}_h$ consisting of
tetrahedral elements.

We discretize the models by the linear finite element methods, and we apply the PCG method with the proposed preconditioner to solve the resulting algebraic systems. The PCG iteration is terminated
in our experiments when the relative residual is less than $10^{-6}$. We will report the iteration counts in the rest of this section.

\subsection{Tests for linear elasticity problems}

In this subsection,we choose the right-hand side $\bm{f}$ of system (\ref{eq:5.1}) such that the analytic solution $\u=(u_1, u_2, u_3)^T$ is given by:
\begin{eqnarray*}
 u_1 &=& x(x-1)y(y-1)z(z-1)\\
 u_2 &=& x(x-1)y(y-1)z(z-1)\\
 u_3 &=& x(x-1)y(y-1)z(z-1)
 \end{eqnarray*}
where the coefficients $\lambda(x) = \mu(x) = 1$. In our experiments, the
right-hand side $\bm{f}$ is fixed.

In this part, we test the action of the preconditioner $B$ described by {\bf
Algorithm 3.1}. We consider the following  different distributions of the coefficients $\lambda(x)$, $\mu(x)$:

\medskip
\noindent{\bf Case (i)}: the coefficients have no jump, i.e., $\lambda(x) = \mu(x) = 1$ .

\medskip
\noindent{\bf Case (ii)}: the coefficients have large jumps, i.e.,
\[\lambda(\x)=\left\{
\begin{array}{l}
\lambda_0,\quad\quad in~D\cr  \\
~~1,~\quad in~\o\backslash D,
\end{array}
\right. \quad\quad  \mu(\x)=\left\{
\begin{array}{l}
\mu_0,\quad\quad in~D,\cr  \\
~~1,~\quad in~\o\backslash D
\end{array}
\right.
\]
Hereafter
we consider two choices of $D$:

$$ \mbox{Choice (1)}. ~~~D=[\dfrac{1}{4},~\frac{1}{2}]^3;~~~~\mbox{Choice (2)}.~~~ D=[\frac{1}{4},~\frac{1}{2}]^3\bigcup[\frac{1}{2},~\frac{3}{4}]^3. $$
The iteration counts of the PCG method with $B$ are listed in \textsc{Table} \ref{E-tp-I-sm} (for {\bf Case} (i)) and \textsc{Table} \ref{E-tp-I-nsm} (for {\bf Case} (ii)).

%

\vskip 0.1in
  \begin{center}  
  \tabcaption{}\label{E-tp-I-sm}
  Iteration counts of PCG with the preconditioner $B$: the coefficients have no jump
  \vskip 0.2in
\begin{tabular}{|c|c|c|c|c|}
\hline
$ m \backslash n$  &\quad 4 \quad\quad & \quad 6 \quad\quad  & \quad 8 \quad\quad & \quad 10 \quad\quad         \\\hline
                 4&   18&       19&        19&          18 \\\hline
                 8&   20&       20&       20&       19     \\\hline
                12&  22&       21&      21&       20 \\\hline
                16&  23&       23&      22&       21 \\\hline
\end{tabular}
\end{center}
 \vskip 0.1in

\vskip 0.1in
\begin{center}
  \tabcaption{}\label{E-tp-I-nsm}
   Iteration counts of PCG with the preconditioner $B$: the coefficients have large jumps
  \vskip 0.2in
\begin{tabular}{|c|c|c|c|c|c|c|c|c|}
\hline
$\mbox{}$
&\multicolumn{4}{c|}{Choice~(1)~of~$D$}
&\multicolumn{4}{c|}{Choice~(2)~of~$D$}\\\hline
$\mbox{}$
&\multicolumn{2}{c|}{$ \lambda_0 =\mu_0=10^{-5}$}
&\multicolumn{2}{c|}{$ \lambda_0 =\mu_0=10^{5} $}
&\multicolumn{2}{c|}{$ \lambda_0 =\mu_0=10^{-5} $}
&\multicolumn{2}{c|}{$ \lambda_0 =\mu_0=10^{5}  $}\\\hline
$~m ~ \backslash ~ n$&4   &8    &4  &8   &4  &8  &4  &8\\\hline
 4 &16      &18     &25   &22   &16      &19   &25  &22 \\\hline
 8 &17    &20   &27   &23   &17    &21   &27&23  \\\hline
 12&19   &21   &28   &24   &18    &23   &28&24 \\\hline
 16&20   &22  &29    &25   &20    &24   &29&26  \\\hline
\end{tabular}
\end{center}
 \vskip 0.1in

We observe from \textsc{Table} \ref{E-tp-I-sm} that the iteration counts of PCG method grows slowly when $m =d/h$ increases but $n=1/d$ is fixed,
and that these counts vary stably when $m$ is fixed but $n$ increases. This show that, when the coefficients is smooth, the condition number of the preconditioned system $B^{-1}A$ should grow
logarithmically with $d/h$ only, not depend on $1/d$. The data in
 \textsc{Table} \ref{E-tp-I-nsm} indicate that, even if the coefficients have large jumps, the iteration counts of PCG still grows slowly. It confirms that the preconditioner $B$ is effective for the system arising from nodal element discretization for linear elasticity problems.


\subsection{Tests for Maxwell's equations}
In this subsection, we consider Maxwell's equations.
For the time-dependent Maxwell's equations, we need to solve the following {\bf curlcurl}-system at each time step (see \cite{Cessenat1998, Hiptmair2002, Monk2003}):
\begin{equation}
\left \{
\begin{array}{rrr}
\c(\al\,\c\,\u)+\beta\u=\f, &in& \q\o,
\\
\u \times \bm{n} = 0, &on& \partial\Omega
\end{array}
\right.
\label{eq:5.2}
\end{equation}
where the coefficients $\al(\x)$ and $\beta(\x)$ are two positive
bounded functions in $\o$. $\n$ is the unit outward normal vector on $\p\o$.

Let $H({\bf curl};\o)$ be the Sobolev space consisting of all square
integrable functions whose {\bf curl}'s are also square integrable
in $\o$, and $ H_0({\bf curl};\o)$ be a subspace of $H({\bf
curl};\o)$ of all functions whose tangential components vanishing
on $\p\o$. In an analogous way, in order to get the weak form of
(\ref{eq:5.2}), just like linear elasticity problems, let $V(\Omega)=H_0({\bf curl})$, and
\[ \A(\bm{u},\bm{v}) = \int_{\Omega}(\alpha~\c~\u \cdot \c~\v
+ \beta~\u \cdot \v)d\x, \]
\[ \langle \bm{F},\bm{v} \rangle = \int_{\Omega}\f\cdot\bm{v}d\x. \]

Let $R(K)$ be a subset of all linear polynomials on the element
$K$ of the form:
$$ R(K)=\Big\{\a+\b\ti\x;~ \a,\b\in \mathds{R}^3, ~\x\in K\Big\} \,. $$
It is well-known that for any $\v\in V_{h}(\o)$, its tangential
components are continuous on all edges of each element in the
triangulation $\mathcal{T}_h$. Moreover, each edge element function
$\v$ in $V_{h}(\o)$ is uniquely determined by its moments on each edge
$e$ of $\mathcal{T}_h$:
\begin{equation}
 \Big\{\la_e(\v)=\int_e\v\cdot\t_e ds; ~ e\in \mathcal{E}_h\Big\},
 \label{eq:5.3}
\end{equation}
where $\t_e$ denotes the unit vector on the edge $e$.

As in the last subsection, we assume that $\O$ can be written as the union of polyhedral
subdomains $D_1$, $\cdots$, $D_{N_0}$ with $N_0$ being a fixed positive integer, such that
$\alpha(x)= \alpha_r$ and $\beta(x)=\beta_r$ for $x\in D_r$, where
every $\alpha_r$ and $\beta_r$ is a positive constant. Let the subdomains $\Omega_k$ satisfy the condition: each $D_r$ is the union of some subdomains in $\{\Omega_k\}$.

Let the right-hand side $\bm{f}$ in the equations
(\ref{eq:5.2}) to be selected such that
the exact solution $\u=(u_1, u_2, u_3)^T$ is given by
\begin{eqnarray*}
 u_1&=&xyz(x-1)(y-1)(z-1)\,, \\
 u_2&=& \sin (\pi x) \sin(\pi y) \sin(\pi z)\,, \\
 u_3&=& (1-e^x)(1-e^{x-1})(1-e^y)(1-e^{y-1})(1-e^z)(1-e^{z-1})\,,
\end{eqnarray*}
where the coefficients $\alpha(\x)$ and $\beta(\x)$ are both constant $1$.
This right-hand side $\bm{f}$ is also fixed in our experiments.

In this part, we investigate the effectiveness of the preconditioner $B$ described by {\bf
Algorithm 3.1}.
We consider the following  different distributions of the coefficients $\alpha(\x)$ and $\beta(\x)$:

\medskip
\noindent{\bf Case (i)}: the coefficients have no jump, i.e., $\alpha(\x)=\beta(\x)=1$.

\medskip
\noindent{\bf Case (ii)}: the coefficients have large jumps, i.e.,
\[\alpha(\x)=\left\{
\begin{array}{l}
\al_0,\quad\quad in~D\cr  \\
~~1,~\quad in~\o\backslash D,
\end{array}
\right. \quad\quad  \beta(\x)=\left\{
\begin{array}{l}
\beta_0,\quad\quad in~D,\cr  \\
~~1,~\quad in~\o\backslash D
\end{array}
\right.
\]

 The iteration counts of the PCG method with $B$ are listed in
\textsc{Table}  \ref{M-tp-I-sm} (for {\bf Case} (i)) and  \textsc{Table} \ref{M-tp-I-nsm} (for {\bf Case} (ii)).

\vskip 0.1in
\begin{center}
\tabcaption{}\label{M-tp-I-sm}
 Iteration counts of PCG with the preconditioner $B$: the coefficients have no jump
\vskip 0.2in
\begin{tabular}{|c|c|c|c|c|c|}
\hline
$ m \backslash n$ &\quad 4\quad\quad       &\quad 6 \quad\quad   &\quad 8 \quad\quad         &\quad 10 \quad\quad \\\hline
  4  &  16      &16       &15      &15 \\\hline
  8  &  18      &17       &17      &16 \\\hline
 12 &  19      &19       &18      &18\\\hline
 16 &  20      &20       &19      &19\\\hline
\end{tabular}
 \end{center}
  \vskip 0.1in

\vskip 0.1in
\begin{center}
\tabcaption{}\label{M-tp-I-nsm}
 Iteration counts of PCG with the preconditioner $B$: the coefficients have large jumps
\vskip 0.2in
\begin{tabular}{|c|c|c|c|c|c|c|c|c|}\hline
$\mbox{}$& \multicolumn{4}{c|}{Choice~(1)~of~$D$}&
\multicolumn{4}{c|}{Choice~(2)~of~$D$}\\\hline $\mbox{}$&
\multicolumn{2}{c|}{$ \beta_0 =\alpha_0 =10^{-5}  $
}&\multicolumn{2}{c|}{$ \beta_0 =\alpha_0=10^{5} $}&
\multicolumn{2}{c|}{$ \beta_0 =\alpha_0=10^{-5}  $
}&\multicolumn{2}{c|}{$ \beta_0 =\alpha_0=10^{5} $}\\\hline
$m ~ \backslash n$& 4  &8  &4  &8   &4   &8   &4   &8\\\hline
 4& 14     &15  &19  &18       &13    &16    &18    &19    \\\hline
 8& 15     &17  &21  &20     &15      &17    &21  &21    \\\hline
12&16     &19  &23  &21     &16     &19     &22  &23    \\\hline
16&17     &20  &24  &22     &17     &20     &24  &24    \\\hline
\end{tabular}

\end{center}
 \vskip 0.1in

We observe from the above two tables that, although the coarse
space is chosen as the simplest one for Maxwell's
equations, the iteration counts of the PCG method with the preconditioner $B$ grow logarihmically with $m=d/h$ only, not depend on $n=1/d$, even if the coefficients have large jumps.

\subsection{Numerical results on irregular subdomains}
In this subsection, we consider the case of irregular subdomains explained in
Subsection 3.4. As usual (refer to \cite{DohrmannW2015}), here we consider only the case with constant coefficient 1
(if the coefficients have large jumps, we can not require that the distribution of the jumps of the coefficients is consistent with that of the irregular subdomains).

We still use $N$ to denote the number of subdomains (which are not cubes any more)
and $h$ to denote the fine mesh size.

Firstly, we consider the domain $\Omega$ is a unit cube (see \textsc{Figure} \ref{fig:cube_dom}).
We do the experiments for both linear elasticity problems and Maxwell's equations.
In order to describe the results clearly, we use the
previous symbols $m$ and $n$.  Here $n = \sqrt[3]{N}$ and $m  =\frac{1}{h\sqrt[3]{N}}$.
In \textsc{Table} \ref{cube-ire}, we list the iteration counts of PCG method with the proposed preconditioner.

\begin{center}
\tabcaption{}\label{cube-ire}
 Iterations counts of PCG method with the preconditioner $B$ for irregular subdomain partition on cubic domain
 \vskip 0.2in
\begin{tabular}{|c|c|c|c|c|}
\hline
\multicolumn{5}{|c|}{\mbox Maxwell's equations \mbox} \\\hline

$ m\backslash n $ &   $\quad 4 \quad$ & $\quad 6\quad$ & $\quad 8\quad$ & $\quad 10 \quad$ \\\hline
$4$& 21  & 21 & 21   & 22 \\\hline
$8$& 23  & 24 & 22   & 24 \\\hline
$12$&25 &25 &25     &25   \\\hline
$16$&27 &26 & 27   & 26 \\\hline
\multicolumn{5}{|c|}{\mbox Linear elasticity problems \mbox}\\\hline
$ m\backslash n $ &   $4$ & $6$&  $8$ & $10$ \\\hline

 $4$& 24 & 23 & 24 & 23 \\\hline
 $8$&28 & 29 & 27 & 29 \\\hline
 $12$&30 & 32 &33   & 32 \\\hline
 $16$&32& 33 &33 & 34 \\\hline

\end{tabular}

\end{center}
 \vskip 0.1in

%
%
%

\begin{figure}[hbt]
\centering
\begin{tabular}{c}
\includegraphics[width =8cm]{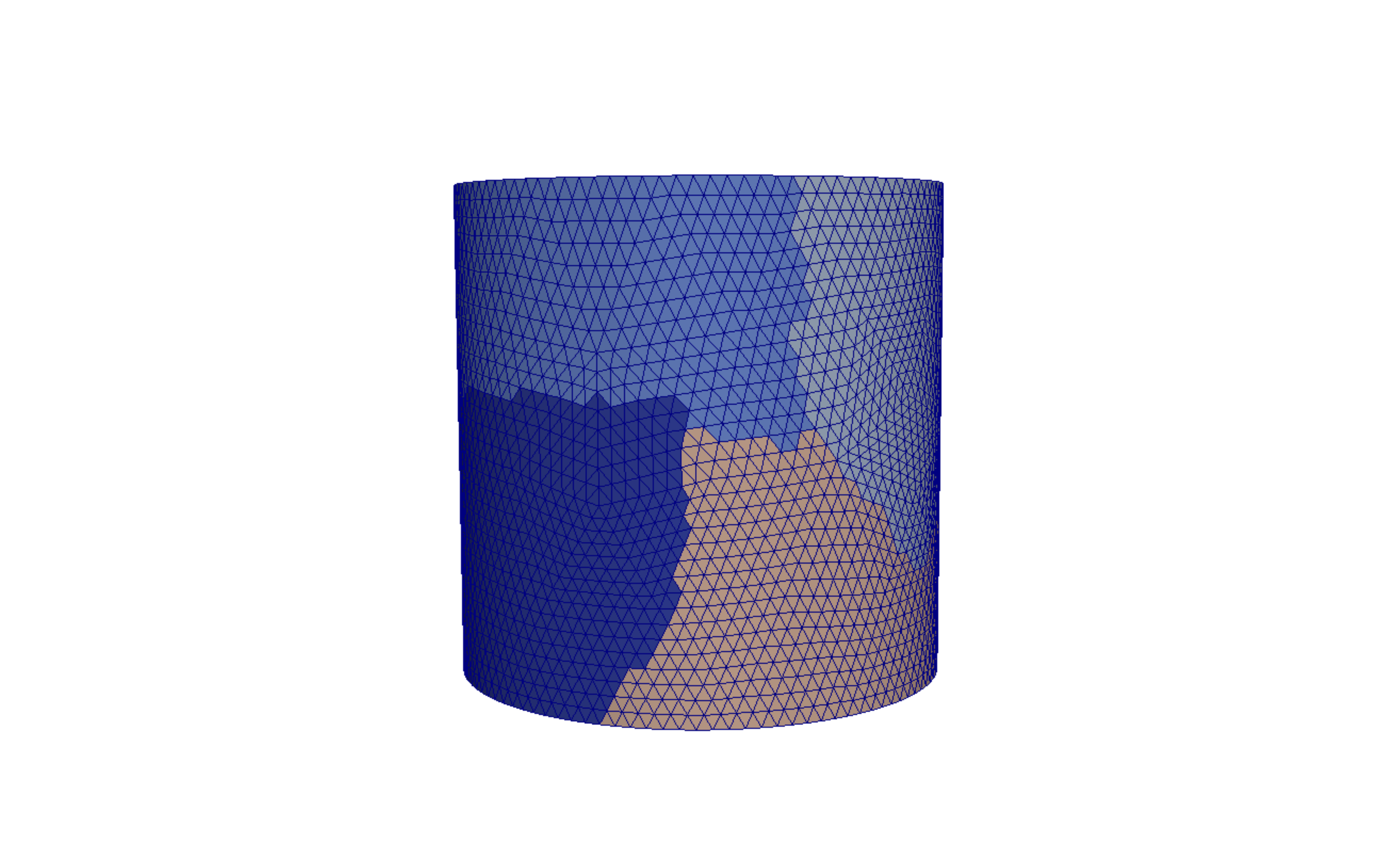}
\end{tabular}
\caption{irregular partition for $\Omega$}
\label{fig:cylinder_dom}
\end{figure}

Secondly, we consider the standard cylinder domain. The mesh partition and subdomains are shown in \textsc{Figure} \ref{fig:cylinder_dom}. We use $N_e$ to denote the
total number of elements of fine mesh. We test two series of mesh: $N_e =765952$ and $N_e = 6127616$.  Here $n = \sqrt[3]{N}$ and and $m \approx \sqrt[3]{{N_e\over N}}$.
The iteration counts of PCG are listed in  \textsc{Table} \ref{cylinder-ire}.

%
%
%

\begin{center}
\tabcaption{}\label{cylinder-ire}
 Iterations counts of PCG method with the preconditioner $B$ for irregular subdomain partition on cylindrical domain
 \vskip 0.2in
\begin{tabular}{|c|c|c|c|c|c|c|c|}
\hline
\multicolumn{8}{|c|}{Maxwell's equations} \\\hline
$ m \backslash n $ & 4   & $ m \backslash n $&  6  & $ m \backslash n $& 8  &  $ m \backslash n $  &10 \\\hline
 $23$   &24   &$15$     &23   &$11$ &22   &$9$   &22  \\\hline
 $46$   &23   &$31$     &23   &$23$ &26   &$18$ &23\\\hline
\multicolumn{8}{|c|}{Linear elasticity problem}\\\hline
$ m \backslash n $ & 4   & $ m \backslash n $&  6  & $ m \backslash n $& 8  &  $ m \backslash n $  &10 \\\hline
 $23$   &28   &$15$     &26   &$11$ &24   &$9$   &25  \\\hline
 $46$   &39   &$31$     &32   &$23$ &29   &$18$ &27 \\\hline
\end{tabular}

\end{center}
 \vskip 0.1in

The above two tables indicate that, even if the subdomains are irregular,  the proposed preconditioner is still effective.
We point out that the dimensions of the coarse space in the current situation have not obvious increase comparing with the case with regular subdomains considered
in previous two subsections.

\subsection{Comparison with the  BDDC method}
In this part, we compare the proposed method with the BDDC method.  All our codes run in sequential way, except when coarse problem and local problems are solved by direct solver package MUMPS, in which
 high performance OpenBLAS libraries based on multithreading parallelism are called.
  In \textsc{Table} \ref{E-comp}, we list the PCG iteration counts and the
 running time for the linear elasticity problems and in \textsc{Table} \ref{M-comp}, \ref{M-irre-comp} for Maxwell's equations.
 Here we use $B$ and $\tilde{B}$ to denote the preconditioners with exact vertex-related solvers $B_{\vv}$ and with inexact one $\tilde{B}_{\vv}$, respectively.

 In order to shorten the length of the paper, in \textsc{Table} \ref{E-comp} we only list the results for $\lambda_0 = \mu_0$ and regular subdomain partition. The parameters $m$ and $n$ are the same meaning as in \textsc{Table} \ref{E-tp-I-nsm}.
 As mentioned in \cite{Dohrmann2003},  when using corner constraints only, we get bad performance in the BDDC method. Because of this, in \textsc{Table} \ref{E-comp} we consider both corner and edge constraints for
 the BDDC method.

 \vskip 0.1in
\begin{center}
  \tabcaption{}\label{E-comp}
   Iteration counts and time(s) for the linear elasticity  problems:  Choice (2) of $D$,  $\lambda_0 = \mu_0$ and  regular partition on cubic domain
  \vskip 0.2in
  \tiny
\begin{tabular}{|c|c|c|c|c|c|c|c|c|c|c|c|c|c|}
\hline
& &\multicolumn{4}{|c|}{BDDC}
 &\multicolumn{4}{|c|}{$B$}
 &\multicolumn{4}{|c|}{$\tilde{B}$} \bigstrut \\\hline
$\lambda_0$& $m=n$ &iter &setup  &solve & total &iter &setup  &solve & total &iter &setup  &solve & total \bigstrut  \\\hline
\multirow{3}[6]*{$10^5$}
& 4  &7 &0.6 &0.1   &0.7     &25 &0.4  &0.2  &0.6       &25  &0.3  &0.2   &0.5  \bigstrut \\\cline{2-14}
& 8  &11 &52.2 &10.5 &62.7     &23 &24.6 &19.8 &44.4      &25 &16.7 &14.4  &31.1 \bigstrut \\\cline{2-14}
&12 &14 &696.6 &179.1 &875.7 &24 &406.8 &290.8 &697.6 &26 &256.0 &221.4 &477.4 \bigstrut  \\\hline
\multirow{3}[6]*{$1$}
 &4  &7 &0.5 &0.1   &0.6     &18 &0.4  &0.2  &0.6       &19  &0.3  &0.1   &0.4  \bigstrut \\\cline{2-14}
& 8  &10 &57.4 &11.3 &68.7     &20 &25.4 &17.2 &42.6      &21 &16.5 &12.2  &28.7 \bigstrut \\\cline{2-14}
&12 &13 &734.1 &167.3 &901.4 &20 &356.5 &244.4 &600.9 &22 &254.0 &189.9 &443.9 \bigstrut  \\\hline
\multirow{3}[6]*{$10^{-5}$}
 &4  &7 &0.6 &0.1   &0.6     &16 &0.4  &0.1  &0.5       &16  &0.3  &0.1   &0.4 \bigstrut  \\\cline{2-14}
& 8  &11 &52.9 &10.6 &63.5     &21 &24.7 &18.1 &42.8      &23 &17.6 &13.3  &30.9 \bigstrut \\\cline{2-14}
&12 &13 &716.1 &167.2 &883.3 &23 &364.7 &286.0 &650.7 &26 &255.2 &219.0 &474.2 \bigstrut  \\\hline
\end{tabular}
\end{center}
\vskip 0.1in

It can be seen from \textsc{Table} \ref{E-comp} that the BDDC preconditioner has the smallest PCG iteration counts among all the tested preconditioners. However, one has to solve algebraic
equations to get coarse basis functions in the BDDC method. As mentioned in the existing works, if only corner constraints are used, the
iteration counts varies unstably for the BDDC method. When edges and corner constraints are used, one needs to solve local saddle-point problem 60 times for each floating
subdomain in the BDDC method.
In addition, the adding constraints to local problems enlarge the size of local problems for BDDC mehod. This can explain why
the total CPU time for setting-up in the BDDC method is more than that in the proposed method. Due to these  reasons, the proposed method has less cost than the BDDC method.

As for Maxwell's equations, we list the results in \textsc{Table} \ref{M-comp} for regular subdomain partition and in \textsc{Table} \ref{M-irre-comp} for irregular subdomain partition. As mentioned in \cite{DohrmannW2015},  instead of using deluxe scaling,  we use e-deluxe scaling in BDDC method, which can result in
significant computional savings and indistinguishable iteration counts.

\vskip 0.1in
\begin{center}
  \tabcaption{}\label{M-comp}
   Iteration counts and time(s) for Maxwell's equation: $\beta=1$ and  regular partition on cubic domain
  \vskip 0.2in
  \tiny
\begin{tabular}{|c|c|c|c|c|c|c|c|c|c|c|c|c|c|}
\hline
 &&\multicolumn{4}{|c|}{BDDC}
 &\multicolumn{4}{|c|}{$B$}
 &\multicolumn{4}{|c|}{$\tilde{B}$} \bigstrut\\\hline
$\alpha $&$n=m$ &iter &setup  &solve & total &iter &setup  &solve & total &iter &setup  &solve & total \bigstrut \\\hline
\multirow{3}[6]*{$10^3$}
  &4  &9 &0.8 &0.4    &1.2  &17 &0.7  &0.4  &1.1   &19  & 0.6  &0.5   &1.1 \bigstrut  \\\cline{2-14}
 & 8  &12 &73.7 &25.3 &99.0  &19 &45.0  &36.6 &81.6 &21 &32.8 &27.4 &59.2 \bigstrut  \\\cline{2-14}
 &10 &13 &337.7 &106.8 &444.5 &19 &193.8 &148.1 &341.9 &21 &137.4 &122.8 &260.2 \bigstrut  \\\hline
\multirow{3}[6]*{$10^2$}
&  4  &9 &0.8 &0.3    &1.1  &17 &0.6  &0.4  &1.0   &19  & 0.7  &0.5   &1.1 \bigstrut  \\\cline{2-14}
&  8  &12 &62.5 &24.5 &87.0  &18 &44.3  &34.3 &78.7 &21 &33.2 &28.4 &61.6  \bigstrut \\\cline{2-14}
& 10 &13 &291.6 &111.7 &403.3 &19 &192.5 &158.8 &351.3 &21 &138.5 &123.8 &262.5 \bigstrut  \\\hline
\multirow{3}[6]*{$1$}
&  4  &9 &0.9 &0.3    &1.2  &16 &0.7  &0.4  &1.1   &18  & 0.7  &0.5   &1.2 \bigstrut \\\cline{2-14}
 & 8  &12 &62.7 &24.9 &87.6  &17 &49.0  &37.1 &86.1 &19 &33.3 &24.5 &57.8 \bigstrut \\\cline{2-14}
& 10 &12 &296.8 &102.5 &399.3 &18 &198.5 &154.1 &352.6 &19 &138.5 &112.2 &250.7 \bigstrut  \\\hline
\multirow{3}[6]*{$10^{-4}$}
&  4  &9 &0.8 &0.3  &1.1  &11 &0.6  &0.3  &0.9   &18  & 0.6  &0.4   &1.0 \bigstrut \\\cline{2-14}
&  8  &12 &64.0 &24.9 &88.8  &11 &47.1  &21.0   &68.10 &17 &32.1 &21.7 &53.8 \bigstrut \\\cline{2-14}
& 10 &13 &296.8 &102.5 &399.3 &11 &198.5 &154.1 &352.6 &19 &138.5 &112.2 &250.7 \bigstrut  \\\hline
\multirow{3}[6]*{$10^{-6}$}
 & 4  &13 &0.8 &0.5  &1.3            &11 &0.7  &0.3  &1.0   &20  & 0.5  &0.5   &1.0  \bigstrut \\\cline{2-14}
 & 8  &19 &75.7 &38.4 &114.1     &9 &46.5  &16.3   &62.8 &22 &32.0 &28.4 &60.4 \bigstrut \\\cline{2-14}
 &10 &20 &274.2 &152.9 &427.1 &9 &192.9 &74.7 &267.6 &17 &140.9 &98.8 &239.7 \bigstrut  \\\hline
\end{tabular}
\end{center}

\vskip 0.1in

In the case of regular subdomain partition, for each floating cubic subdomain, one needs to solve local saddle-point problems 24 times to get coarse basis functions in the BDDC method. Because of this, the results in \textsc{Table} \ref{M-comp} indicate that the setting up time of BDDC method is larger than the proposed method. Then
the total time of BDDC method is larger than the proposed methods although its iteration counts is smaller than our method.  When solving $B_{\vv}$ in inexact manner, we can further
save much time than $B$ in the setting up phase.  From \textsc{Table} \ref{M-comp}, we observed that the preconditioner $\tilde{B}$ has faster convergence than
the BDDC preconditioner. We know that, if $\alpha / \beta << 1$, the operator $A$ has very good conditioning. But, we found a surprising phenomenon: when increasing the value of $m$, the iteration counts of the BDDC method
grows quickly and is larger than the proposed method.

\vskip 0.1in
\begin{center}
  \tabcaption{}\label{M-irre-comp}
   Iteration counts and time(s) for Maxwell's equation: $\beta=1$ and irregular partition on cubic domain
  \vskip 0.2in
\begin{tabular}{|c|c|c|c|c|c|c|c|c|c|}
\hline

& &\multicolumn{4}{|c|}{BDDC}
  &\multicolumn{4}{|c|}{B} \bigstrut \\\hline
$\alpha$& $n=m$ &iter &setup  &solve & total &iter &setup  &solve & total \bigstrut \\\hline
\multirow{3}[6]*{$10^3$}
 & 4  &13 &1.6 &1.4 &3.0            &23 &1.1  &0.8  &1.9 \bigstrut  \\\cline{2-10}
 & 8  &18 &112.6 &75.5 &188.1     &25 &69.7  &51.6 &121.3  \bigstrut  \\\cline{2-10}
 &10 &20 &418.9&252.8 &671.7 &27 &293.9 &184.4 &478.3 \bigstrut   \\\hline
\multirow{3}[6]*{$10^2$}
 & 4  &13 &1.7 &1.5  &3.1             &23 &1.0  &0.8  &1.8 \bigstrut \\\cline{2-10}
 &8  &18 &114.8 &75.4 &190.2    &24 &69.0  &49.5 &118.5 \bigstrut \\\cline{2-10}
 &10 &20 &421.4&249.3 &670.7  &26 &296.0 &177.9 &473.9 \bigstrut \\\hline
\multirow{3}[6]*{$1$}
&  4  &13 &1.6 &1.4  &3.0             &21 &1.0  &0.7  &1.7 \bigstrut \\\cline{2-10}
 & 8  &18 &111.8 &75.2 &187.0    &22 &68.8  &45.9 &114.7 \bigstrut \\\cline{2-10}
 &10 &20 &417.7&254.8 &672.5  &24 &299.5 &165.7 &465.2 \bigstrut \\\hline
\multirow{3}[6]*{$10^{-4}$}
&  4  &11 &1.6 &1.3  &2.9             &13 &1.1  &0.5  &1.6 \bigstrut \\\cline{2-10}
&  8  &15 &113.3 &63.6 &177.0    &13 &69.4  &26.9 &97.3 \bigstrut \\\cline{2-10}
& 10 &19 &422.6&238.6&661.2  &14 &293.5 &97.9 &391.4 \bigstrut \\\hline

\multirow{3}[6]*{$10^{-6}$}
&  4  &12 &1.6 &1.3  &2.9             &14 &1.1  &0.5  &1.5 \bigstrut\\\cline{2-10}
& 8  &21 &113.5 &87.5 &201.0    &12 &69.3  &26.8 &96.1\bigstrut \\\cline{2-10}
& 10 &28 &420.6 &343.7&764.2    &11 &293.9 &77.8 &371.7 \bigstrut \\\hline
\end{tabular}
\end{center}
\vskip 0.1in

In  \textsc{Table} \ref{M-irre-comp}, the parameters $m$ and  $n$ have the same meanings as in \textsc{Table} \ref{cube-ire}. We observed that both setting up time and PCG solving time 
in the proposed method are less than that in the BDDC method. Since the geometric properties of subdomains  generated by Metis are very complicated,
the number of subdomain edges on each irregular subdomain may be much larger than the one in a cubic subdomain and so the total subdomain edges is much larger than the regular subdomain partition case. 
The drawback in irregular subdomain partition enlarges the size of local saddle-point problems and increase the cost of calculation for e-deluxe scaling operator.
This investigation can explain why the PCG solving time in the BDDC method is also larger than that in the proposed method. For the case with irregular subdomains, we have not tested the preconditioner 
$\tilde{B}$ with coarsening technique since the existing coarsening software is not practical yet. In near future, we will try to design special coarsening software for irregular subdomains.

\section{Conclusion}
 In this paper, we have
constructed one substructuring preconditioner with the simplest coarse space for
general elliptic-type problems in three dimensions. In particular, we design 
new local interface solvers, which are easy to implement and do not depend on the considered models.
The proposed preconditioner can absorb some advantages of the non-overlapping
and overlapping domain decomposition methods. We have given an analysis of convergence of the preconditionner for
linear elasticity problems, which shows the proposed preconditonner is  nearly optimal  and also robust with respect to
the (possibly large) jumps of the coefficient. We also consider the case with irregular subdomain partition for numerical experiments.
 We have given some numerical results to show that the proposed preconditionner is effective uniformly for
the linear elasticity problem and Maxwell's equations in three dimensions.


\bibliographystyle{siamplain}

\end{document}